\documentclass[georeport,12pt]{geophysics}

\usepackage{amsmath}
\usepackage{amssymb}
\usepackage{amsthm}
\usepackage{amscd}
\usepackage{setspace}
\usepackage{wasysym}
\usepackage{algorithm}
\usepackage{algpseudocode}

\newcommand{\bx}{\mathbf{x}}

\newcommand{\bv}{\mathbf{v}}
\newcommand{\bff}{\mathbf{f}}
\newcommand{\bu}{\mathbf{u}}

\newcommand{\bR}{\mathbf{R}}

\setfigdir{.}

\begin{document}
\title{Efficient Computation of Extended Surface Sources}
\author{William W. Symes\\
  Department of Computational and Applied Mathematics,\\
  Rice University,Houston TX 77251-1892 USA,\\
  email {\tt symes@rice.edu},\\
ORCID 0000-0001-6213-4272}

\lefthead{Symes}

\righthead{Approximate Source Inversion}

\maketitle
\begin{abstract}
Source extension is a reformulation of inverse problems in wave propagation, that at least in some cases leads to computationally tractable iterative solution methods. The core subproblem in all source extension methods is the solution of a linear inverse problem for a source (right hand side in a system of wave equations) through minimization of data error in the least squares sense with soft imposition of physical constraints on the source via an additive quadratic penalty. A variant of the time reversal method from photoacoustic tomography provides an approximate solution that can be used to precondition Krylov space iteration for rapid convergence to the solution of this subproblem. An acoustic 2D example for sources supported on a surface, with a soft contraint enforcing point support, illustrates the effectiveness of this preconditioner.
\end{abstract}


\inputdir{.}

\section{Introduction}
Full Waveform Inversion (FWI) can be described in terms of 
\begin{enumerate}
\item a linear wave operator $L[{\bf c}]$, depending on a vector of
  space-dependent coefficients ${\bf c}$ and acting on causal vector wavefields $\bu$ vanishing in negative time:
\begin{equation}
\label{eqn:init}
\bu \equiv 0, t \ll 0; 
\end{equation}
\item a trace sampling operator $P$ acting on wavefields and producing data traces;
\item and a (vector) source function (of space and time) $\bff$ representing energy input to the system. 
\end{enumerate}
The basic FWI problem is: given data $d$, find ${\bf c}$ so that 
\begin{equation}
\label{eqn:fwi}
P\bu \approx d \mbox{ and } L[\bf{c}]\bu = \bff.
\end{equation}
In this formulation, the source function $\bff$ may be given, or
to be determined subject to some constraints.

A simple nonlinear least squares formulation is:
\begin{equation}
\label{eqn:ols}
\mbox{choose } {\bf c} \mbox{ to minimize } \|PL[{\bf c}]^{-1}\bff -d \|^2.
\end{equation}
Practical optimization formulations typically augment the objective in
\ref{eqn:ols} by additive penalties or other constraints.

As is well-known, local optimization methods are the only feasible
approach given the dimensions of a typical instance of \ref{eqn:fwi},
and those have a tendency to stall due to ``cycle-skipping''. See for
example \cite{VirieuxOperto:09} and many references cited there. Source
extension is one approach to avoiding this problem. It consists in
imposing the wave equation as a soft as opposed to hard constraint, by
allowing the source field $\bff$ to have more degrees of freedom than
is permitted by a faithful model model of the seismic experiment, and
constraining these additional degrees of freedom by means of an
additive quadratic penalty modifying the problem \ref{eqn:ols}:
\begin{equation}
\label{eqn:esi}
\mbox{choose } {\bf c}, \bff \mbox{ to minimize } \|PL[{\bf c}]^{-1}\bff -d \|^2 + \alpha^2 \|A\bff\|^2 
\end{equation}
The operator $A$ penalizes deviation from known (or assumed)
characteristics of the source function - its null space consists of
feasible (or ``physical'') source models.

\cite{HuangNammourSymesDollizal:SEG19} present an overview of the
literature on source extension methods, describing a variety of
methods to add degrees of freedom to physical source model. The present paper
concerns {\em surface source extension}: physical sources are
presumed to be concentrated at points $\bx_s$ in space, whereas their extended
counterparts are permitted to spread energy over surfaces containing
the physical source locations. A simple choice for the penalty
operator $A$ is then multiplication by the distance $|\bx-\bx_s|$ to the physical
source location:
\begin{equation}
  \label{eqn:penop}
  (A\bff)(\bx,t) = |\bx-\bx_s|\bff(\bx,t)
\end{equation}
I shall use this choice of penalty operator whenever a specific choice
is necessary in the development of the theory below.

This paper presents a numerically efficient approach to solving the
{\em source subproblem} of problem \ref{eqn:esi}:
\begin{equation}
\label{eqn:esis}
\mbox{given } {\bf c}, \mbox{ choose } \bff \mbox{ to minimize }
\|PL[{\bf c}]^{-1}\bff -d \|^2 + \alpha \|A\bff\|^2 
\end{equation}
Solution of this subproblem is an essential component of {\em variable
  projection} algorithms for solution of the nonlinear inverse problem
\ref{eqn:esi}. Variable projection is not merely a convenient choice
of algorithm for this purpose: it is in some sense essential, see for
example \cite{Symes:SEG20}. It replaces the nonlinear
least squares problem \ref{eqn:esi} with a {\em reduced} problem, to
be solved iteratively. Each iteration involves solution of the
subproblem \ref{eqn:esis}. Therefore efficient solution of the
subproblem is essential to efficient solution of the nonlinear problem
via variable projection.

The penalty operator $A$ defined in \ref{eqn:penop} is linear, so the source
subproblem is a linear least squares problem. Under some additional
assumptions to be described below, I shall show how to construct an
accurate approximate solution operator for problem
\ref{eqn:esis}. This approximate solution operator may be used to
accelerate Krylov space methods for the solution of the surface source
subproblem \ref{eqn:esis}. Numerical examples suggest the
effectiveness of this acceleration.

I will fully describe a preconditioner for a special
case of the source subproblem \ref{eqn:esis}, in which ${\bf u}$ is an
acoustic field, $L[{\bf c}]$ is the wave operator of linear
acoustodynamics, the spatial positions of traces extracted by $P$ lie
on a depth plane $z=z_r$, and the positions at which the extended
source $\bff$ is nonzero lie on another, parallel, depth plane $z=z_s$. 
This
``crosswell'' configuration simplifies the analysis underlying the
construction of approximate solutions for the source subproblem
\ref{eqn:esis}. It is only one of many transmission configurations for
which similar developments are possible. Perhaps the most important
alternative example is the diving wave configuration, which plays a
central role in contemporary FWI. 

The preconditioner construction is very similar to the time-reversal
method in photoacoustic tomography
\cite[]{StefanovUhlmannIP:09}. Preconditioning amounts to a change of
norm in the domain and range spaces of the modeling operator. In this
case, the modfied norms are weighted $L^2$, and the weight operators
map pressure to corresponding surface source on the source and
receiver planes. This pressure-to-source map is closely related to the ``hyperbolic
Dirichlet-to-Neumann'' operator that plays a prominent role in
photoacoustic tomography and other wave inverse problems
\cite[]{Rachele:00,StefUhl:05}. \cite{HouSymes:EAGE16} demonstrated a
very similar preconditioner for Least Squares Migration, also for its
subsurface offset extension \cite[]{HouSymes:16}, motivated by
\cite{tenKroode:12}. These constructions also involve the
Dirichlet-to-Neumann operator. This concept also turns up in hidden
form in the work of Yu Zhang and collaborators on true amplitude
migration
\cite[]{YuZhang:14,TangXuZhang:13,XuWang:2012,XuZhangTang:11,Zhang:SEG09}.

The obvious computation of the pressure-to-source map - prescribe the
pressure, solve the wave equation with this boundary condition, read
off the equivalent source - suffers from intrinsic numerical
inaccuracy. I suggest an alternative computationally feasible
approach, via economical short-distance wave propagation. Since the
map is symmetric only in an approximate, asymptotic sense, it must be
symmetrized for use as a Krylov preconditioner. I describe a
symmetrization procedure that requires no further wave computations
beyond those necessary to compute the action of the operator itself.

The discussion in this paper is formal and incomplete, in the sense
that some important mathematical underpinnings are only
sketched. I will treat the modeling operator $PL[{\bf c}]^{-1}$ as if
it mapped square integrable surface sources to square integrable sampled
data. This is not true in full generality: while the surface source
problem has distribution solutions, they are not generally square
integrable (finite acoustic field energy). Even if the solutions have
finite energy, they do not in general have well-defined
restrictions to lower-dimensional sets. In other words, the action of
the sampling
operator $P$ on the space-time plane $z=z_r$ is not well-defined for arbitrary
finite-energy acoustic fields. Thus the modeling operator envisioned above
may not be well-defined.

This
phenomenon is related to the ill-posedness of wave equations as
evolution equations in spatial variables, an observation attributed to
Hadamard (see \cite{CourHil:62}, Chapter 6, section 17). Some constraint on the acoustic field,
beyond finite energy, is mandatory in any precise mathematical formulation the inverse problems
\ref{eqn:esi} and \ref{eqn:esis}. In fact, the natural constraint in
the ``crosswell'' geometry of this paper is that high-frequency energy
travel {\em only} along rays crossing the
surfaces $z=z_s, z=z_r$ transversally. That is, source functions on
$z=z_s$ generate waves with energy traveling along rays leaving the surface at a non-zero angle, and
energy arrives at the recording surface $z=z_r$ along rays making
non-zero angle with it. I will call sources, sampled data, and
acoustic fields with this property {\em
  downgoing} (even though the concept also encompasses {\em upcoming}
propagation). Note that the downgoing property restricts the behaviour of
acoustic fields near the source and receiver surfaces - what the
fields do elsewhere is their own business.

Several works have explored the mathematics of the
downgoing condition and its consequences in the context of the scalar second order wave equation, see for instance
\cite{Payn:75,Symes:83,Lasi:86,LasLionsTrig:86,Lasi:87, BaoSy:91b}. 
Elaboration of these mathematical details is beyond the scope of this
paper, which aims instead to explore the algorithmic
consequences of the mathematical structure implied by the downgoing
condition.

The next section defines the modeling operator $PL[{\bf c}]$, its
adjoint, and important specializations (pressure vs. normal velocity
sources and data). The sections to follow define the
source-to-pressure operator, construct an approximate inverse of the
modeling operator by time reversal (as suggested by work in
photoacoustic tomography), use the source-to-pressure operator to
express the approximate inverse as the modeling operator adjoint in
weighted norms (thus establishing that the modeling operator is {\em approximately
  unitary} in the sense of these norms), explain how to use this
construction to precondition Conjugate Gradient iteration, and
organize the preconditioning computation so as to involve only one
extra and relatively inexpensive
wave propagation calculation. The penultimate section displays simple
2D numerical examples of all of the key steps, culminating in a
comparison of straight vs. preconditioned CG iteration. The paper ends
with a brief discussion-and-conclusion section, reviewing what has
been accomplished and listing a few of the many questions left open.

\section{Operators}

For acoustic wave physics, the coefficient vector is
$\bf{c}=(\kappa,\rho)^T$, with components bulk modulus $\kappa$ and
density $\rho$, and the state vector $\bu=(p,\bv)^T$ consists of
pressure $p$ (a scalar space-time field) and particle velocity $\bv$
(a vector space-time field). The wave operator $L[\bf{c}]$ is:
\begin{equation}
\label{eqn:aweop}
L[\bf{c}]\bf{u} = 
\left(
\begin{array}{c}
\frac{1}{\kappa}\frac{\partial p}{\partial t}  + \nabla \cdot \bv, \\
\rho\frac{\partial \bv}{\partial t} + \nabla p.
\end{array}
\right) 
\end{equation}
That is,
\begin{equation}
  \label{eqn:awemat}
  L[{\bf c}] = \left(
    \begin{array}{cc}
      \frac{1}{\kappa}\frac{\partial}{\partial t} & \nabla \cdot \\
      \nabla & \rho \frac{\partial}{\partial t}
    \end{array}
  \right)
\end{equation}
$L[{\bf c}]$ has a well-defined inverse in the sense of distributions
if it is restricted to either causal or anti-causal vector wavefields.

Most of what follows is valid for any space dimension $n >0$. The
coefficient vector $\bf{c}=(\kappa,\rho)$ is defined throughout space
$\bR^n$, the state vector $\bu$ throughout space-time
$\bR^{n+1}$. Whenever convenient for mathematical manipulations,
$n=3$: for instance, I will write $\bx=(x,y,z)^T$ for the spatial
coordinate vector, and refer to the third (vertical) coordinate of
particle velocity as $v_z$. Examples
later in this paper will use $n=2$ for computational convenience.

Since all of the operators in the discussion that follows depend on
the coefficient vector 
$\bf{c}$, I will suppress it from the notation, for example, $L=L[\bf{c}]$.

The surface source extension replaces point sources on or near a
surface in $\bR^3$ with source functions confined to the same
surface. The simplest example of this extended geometry specifies a
plane $\{(x,y,z,t): z=z_s\}$ at source depth $z_s$ as the surface. For
acoustic modeling, surface sources are combinations of constitutive law
defects and loads normal to the surface, localized on $z=z_s$. That
is, right-hand sides in the system $L\bu=\bff$ take the form
$\bff(\bx,t) = (h_s(x,y,t)\delta(z-z_s),
f_s(x,y,t)\bf{e}_z\delta(z-z_s))^T$ for scalar defect $h_s$ and normal
force $f_s$ ($\bf{e}_z=(0,0,1)$). With the choice $L$ given in
\ref{eqn:awemat}, the causal/anti-causal wave system $L\bu^{\pm}=\bff$
takes the form
\begin{eqnarray}
\label{eqn:awepm}
\frac{1}{\kappa}\frac{\partial p^{\pm}}{\partial t} & = & - \nabla \cdot \bv^{\pm} +
h_s \delta(z-z_s), \nonumber \\
\rho\frac{\partial \bv^{\pm}}{\partial t} & = & - \nabla p^{\pm} +
                                                f_s{\bf e} \delta(z-z_s),\nonumber \\
p^{\pm} & =& 0 \mbox{ for } \pm t \ll 0,\nonumber\\ 
\bv^{\pm} & = & 0 \mbox{ for } \pm t \ll 0.
\end{eqnarray}

\noindent {\bf Remark:} In system \ref{eqn:awepm} and many similar
systems to follow, I will use the shorthand
\[
  p^+ = 0 \mbox{ for } t \ll 0 
\]
to mean that $p^+$ is {\em causal}, that is,
\[
  \mbox{For some } T \in \bR, p^+(\cdot,t) = 0 \mbox{ for all } t <
  T.
\]
Similarly,
\[
  p^- = 0 \mbox{ for } t \gg 0 
\]
signifies that $p^-$ is anti-causal.

Extended forward modeling consists in solving \ref{eqn:awepm} and
sampling the solution components at receiver locations. For
simplicity, throughout this paper I will assume that the receivers are
located on another spatial hyperplane $\{(x,y,z,t): z=z_r\}$ at
receiver depth $z_r>z_s$. The constructions to follow involve interchange
of the roles of $z_s$ and $z_r$ (that is, locating sources on $z=z_r$
and receivers at $z=z_s$), so rather than the sampling operator $P$ of
the introduction, I will denote by $P_s,P_r$ the sampling 
operators on $z=z_s$, $z=z_r$ respectively. In practice, sampling
occurs at a discrete array of points (trace locations) on these
surfaces, and over a zone of finite extent. In this theoretical
discussion, I will neglect both finite sample rate and extent, and
regard the data, for example $P_rp^+$, as continuously sampled and
extending over the entire plane $z=z_r$.

As explained in the Introduction, the downgoing constraint on the square-integrable
source functions $h_s, f_s$ is
essential both for finite energy solutions of the system
\ref{eqn:awepm} to exist, and for these solutions to have well-defined
traces on the receiver surface $z=z_r$. This constraint will be
assumed throughout, often tacitly.
For downgoing solutions of system \ref{eqn:awepm}, the key
components ($p^{\pm}$ and $v^{\pm}_z$) are continuous functions of $z$
in the open slab $z_s<z<z_r$ with well-defined limits at the boundary
planes, but may be discontinuous at the source plane
$z=z_s$. Similarly, the roles of $z_s$ and $z_r$ will be interchanged
in some of the constructions to come, and the corresponding solutions
may be discontinuous at $z=z_r$. Accordingly, interpret $P_s$, $P_r$
as the limit from right and left respectively: for $u=p^{\pm}$ or
$v^{\pm}_z$,
\begin{eqnarray}
  \label{eqn:defsamp}
  P_su(x,y,t) &=& \lim_{z \rightarrow z_s^+} u(x,y,z,t),\nonumber \\
  P_ru(x,y,t) &=& \lim_{z \rightarrow z_r^-} u(x,y,z,t).                  
\end{eqnarray}

The causal/anti-causal vector
modeling operators ${\cal S}^{\pm}_{z_s,z_r}$ are defined in terms of
the solutions $(p^{\pm},\bv^{\pm})$ of the systems \ref{eqn:awepm} by
\begin{equation}
  {\cal S}^{\pm}_{z_s,z_r}(h_s,f_s)^T  = (P_rp^{\pm},P_r v_z^{\pm})^T,
  \label{eqn:fwd}
\end{equation}
The subscript signifies that sources are located on $z=z_s$, the
receivers on $z=z_r$. It is necessary to include this information in
the notation, as versions of ${\cal S}^{\pm}$ with sources and receivers in
several locations will be needed in the discussion below.

\noindent {\bf Remark:} To connect with the formulation presented in
the introduction, note that for continuous $u$,
$P_su(x,y,t)=u(x,y,z_s,t)$, and therefore the adjoint of $P_s$ (in the
sense of distributions) is $P_s^Th(x,y,z,t) =
h(x,y,t)\delta(z-z_s)$. Write ${\cal P}_s = \mbox{diag }(P_s,P_s)$ and
similarly for ${\cal P}_r$. Then
\[
  {\cal S}^{+}_{z_s,z_r} = {\cal P}_r L^{-1}({\cal P}_s)^T,
\]
in which $L^{-1}$ is interpreted in the causal sense, and similarly
for ${\cal S}^{-}$. Sources confined to $z=z_s$ are precisely those
functions (distributions, really) output by ${\cal P}_s^T$, so the
problem statements \ref{eqn:esi} and \ref{eqn:esis} can be rewritten
in terms of ${\cal S}^+_{z_s,z_r}$, with $P$ identified with ${\cal P}_r$.

${\cal S}^{\pm}$ is not stably invertible: its columns are
approximately linearly dependent, as will be verified below. The
diagonal components of ${\cal S}^{\pm}$ thus carry essentially all of
its information, and it is in terms of these that a sensible inverse problem
is defined.

Denote by $\Pi_i, i=0,1$ the projection on the first,
respectively second, component of a vector in $\bR^2$. The 
forward modeling operator from pressure source to pressure trace is
\begin{equation}
  \label{eqn:sdef}
  S^{\pm}_{z_s,z_r} = \Pi_0 {\cal S}^{\pm}_{z_s,z_r} \Pi_0^T 
\end{equation}
and the forward modeling operator from velocity source (normal force)
to velocity trace is
\begin{equation}
  \label{eqn:vdef}
  V^{\pm}_{z_s,z_r} = \Pi_1 {\cal S}^{\pm}_{z_s,z_r} \Pi_1^T 
\end{equation}

With these conventions, we can write the version of the source
subproblem \ref{eqn:esis} studied in this paper as
\begin{equation}
  \label{eqn:esisp}
  \mbox{find }h_s\mbox{ to minimize }\|S^{+}_{z_s,z_r}h_s- d\|^2 +
  \alpha^2\|Ah_s\|^2.
\end{equation}

It follows from the adjoint state method (see Appendix A for details) that
\begin{equation}
  \label{eqn:sadj1}
  ({\cal S}^{\pm}_{z_s, z_r})^T = -{\cal S}^{\mp}_{z_r,z_s}
\end{equation}

Define $R$ to be the {\em time-reversal operator} on functions of
space-time, $Rf(\bx,t) = f(\bx,-t)$, and ${\cal R}$ to be the {\em
  acoustic field time-reversal operator}
\begin{equation}
  \label{eqn:trdef}
  {\cal R} \left(
    \begin{array}{c}
      p\\
      \bv
    \end{array}
  \right) =
  \left(
    \begin{array}{c}
      Rp\\
      -R\bv
    \end{array}
  \right)
\end{equation}
Then 
\begin{equation}
  \label{eqn:trsadj}
  {\cal R}{\cal S}^{\mp} = -{\cal S}^{\pm}_{z_r,z_s}{\cal R}
\end{equation}
Since $R^2 = I$ and ${\cal R}^2 = I$, the identities \ref{eqn:sadj1} and \ref{eqn:trsadj} imply that
\begin{equation} 
  \label{eqn:trtr}
 ({\cal S}^{\pm}_{z_s,z_r})^T = {\cal R}{\cal S}_{z_r,z_s}^{\pm}{\cal R}=
 -{\cal S}^{\mp}_{z_r,z_s}.
\end{equation}
The relation \ref{eqn:trtr} implies that
\begin{eqnarray}
  (S^{\pm}_{z_s,z_r})^T &=& -S^{\mp}_{z_r,z_s} \nonumber\\
                        &=& R S^{\pm}_{z_r,z_s}R, \nonumber\\
    (V^{\pm}_{z_s,z_r})^T &=& -V^{\mp}_{z_r,z_s} \nonumber\\
                        &=& R V^{\pm}_{z_r,z_s}R.
                            \label{eqn:trtrcomp}
\end{eqnarray}

\section{Pressure-to-Source}

Since the system \ref{eqn:awepm} has a unique solution by standard
theory \cite[]{Lax:PDENotes}, the source vector field $(h_s,f_s)$
determines the acoustic field $(p^{\pm},\bv^{\pm})$ in space time, and
in particular the limits from the right at $z=z_s$, $P_sp^{\pm}$ and
$P_sv_z^{\pm}$. This relation is not invertible: it is not possible to
prescribe both pressure and normal velocity on a surface such as
$z=z_s$. So the columns of the matrix operator
${\cal S}^{\pm}_{z_s,z_r}$ must satisfy a linear relation. In this
section I will explain this relation; it involves the {\em
  pressure-to-source} map. This operator also turns out to be the
principal component of a preconditioning strategy for iterative
solution of the optimization problem \ref{eqn:esis}, so I will devote
some effort to its proper definition. It is closely related to the
Dirichlet-to-Neumann operator mentioned in the introduction.

While it is not possible to prescribe both pressure and velocity on
$z=z_s$ in solutions of \ref{eqn:awepm}, it is possible to
prescribe pressure only, for instance: if the function $\phi$ on
the surface $z=z_s$ satisfies suitable conditions, for
example the downgoing constraint mentioned earlier, a unique solution
exists for the acoustic system in both half-spaces $\pm z > z_s$:
\begin{eqnarray}
\label{eqn:awe0}
  \frac{1}{\kappa}\frac{\partial p_{\pm}}{\partial t} & = & - \nabla \cdot \bv_{\pm}, \nonumber \\
  \rho\frac{\partial \bv_{\pm}}{\partial t} & = & - \nabla
                                                    p_{\pm}, \nonumber \\
  p_{\pm} & =& 0,  \mbox{ for } t \ll 0, \nonumber\\ 
  \bv_{\pm} & = & 0 \mbox{ for } t \ll 0, \nonumber\\
  \lim_{z \rightarrow z_s^{\pm}}p_{\pm}(x,y,t,z)& =& \phi(x,y,t).
\end{eqnarray}
Note that the subscript $\pm$ here refers to the sign of $z-z_s$, as opposed
to the superscript ${\pm}$, which refers to the sign of $t$ throughout
this paper.

From the boundary condition (last equation in \ref{eqn:awe0}), one
sees that the pressures $p_{\pm}$ in the two half-spaces have the same
limit at the boundary $z=z_s$. Stick the two half-space
solutions together to form an acoustic field $(p^+,\bv^+)$ in all of
space-time, that is,
\begin{equation}
  \label{eqn:awealt}
  p^+(x,y,z,t) =
  \left\{
    \begin{array}{c}
      p_+(x,y,z,t) \mbox{ if } z>0,\\
      p_-(x,y,z,t) \mbox{ if } z<0,
    \end{array}
  \right.
\end{equation}
and a similar definition for $\bv^+$. Then $p^+$ is continuous across
$z=z_s$, and the boundary condition in system \ref{eqn:awe0} may be
written as $P_sp^+=\phi$.

The same construction can be carried out in the anti-causal sense,
with anti-causal half-space solutions glued together to form a
full-space distribution solution $(p^-,\bv^-)$, with the property that
$p^-$ is continuous across $z=z_s$ and $P_sp^-=\phi$.

The reader may object that the notation $(p^\pm,\bv^\pm)$ is already in
use, for the solution of \ref{eqn:awepm}. This objection is
valid. However, {\em in the sense
  of distributions}, $(p^{\pm},\bv^{\pm})$ as defined in display
\ref{eqn:awealt}, is {\em exactly} the causal solution of \ref{eqn:awepm}
for the choice $h_s = -[v^{\pm}_{z}]|_{z=z_s}, f_s=0$, as follows from a
simple integration-by-parts calculation. So the notation is consistent!

The negative jump $-[v^{\pm}_{z}]|_{z=z_s}$ is thus a function of $\phi$. Define
the {\em pressure-to-source} operator $\Lambda^{\pm}_{z_s}$ by
\begin{equation}
  \label{eqn:deflam}
  \Lambda^{\pm}_{z_s}\phi = -[v^{\pm}_{z}]|_{z=z_s}
\end{equation}
The conclusion: if $h_s = \Lambda^{\pm}_{z_s}\phi$ and $f_s=0$ in the
system \ref{eqn:awepm}, then $\phi=P_sp^{\pm}$.

Otherwise put, $S^{\pm}_{z_s,z_s}\Lambda^{\pm}_{z_s} \phi = \phi$, so
$\Lambda^{\pm}_{z_s}$ is inverse to $S^{\pm}_{z_s,z_s}$. The relation
\ref{eqn:trtrcomp} implies in turn that
\begin{equation}
  \label{eqn:lamadj}
  (\Lambda^{\pm}_{z_s})^T = - \Lambda^{\mp}_{z_s}
\end{equation}

There is also a {\em velocity-to-source} operator. For the solution
$(p^{\pm},\bv^{\pm})$ of system \ref{eqn:awepm} with $h_s=0$, the
normal component of velocity, $v^{\pm}_z$, is continuous across
$z=z_s$, and the velocity source (vertical load)
$f_s=-[p^{\pm}]_{z=z_s}$. I will not name the velocity-to-source
operator, as it does not appear explicitly in the developments to
follow. As will be seen, it is essentially the inverse of the
pressure-to-source operator.

The quadratic form defined by $\Lambda^{\pm}_{z_s}$ has fundamental
physical significance. Define the total acoustic energy $E^{\pm}(t)$ of the
field $(p^{\pm},\bv^{\pm})$, at time $t$ by
\begin{equation}
  \label{eqn:defae0}
  E^{\pm}(t) = \frac{1}{2} \int \,d\bx \, \left(\frac{(p^{\pm})^2}{\kappa} + \rho |\bv^{\pm}|^2\right)(\bx,t).
\end{equation}
Then
\begin{equation}
  \label{eqn:elim}
  \pm \lim_{\pm t \rightarrow \infty} E^{\pm}(t) =  \langle P_sp^{\pm},
  (\Lambda^{\pm}_{z_s} P_sp^{\pm}) \rangle_{L^2(z=z_s)}.
\end{equation}
That is, the value of the quadratic form defined by
$\Lambda^{\pm}_{z_s}$, evaluated at the pressure trace on $z=z_s$,
gives the total energy transferred from the source to the
acoustic field over time. Since $E$ is itself a positive definite
quadratic form in the acoustic field, it follows that $\pm
\Lambda^{\pm}_{z_s}$ is positive semi-definite. 

While $\Lambda^{\pm}_{z_s}$ is positive semi-definite, it is not
symmetric. However, it is {\em approximately symmetric} in the
high-frequency sense. This fact follows from a
geometric optics analysis of the half-space solution. This leads to
the identification of $\Lambda^{\pm}_{z_s}$ as a {\em
  pseudodifferential operator} of order zero on $z=z_s$, with principal symbol
\begin{equation}
  \label{eqn:lamsym}
  \sigma_0(\Lambda^{\pm}_{z_s}) = \pm 2(\kappa(\bx)
  \rho(\bx))^{1/2}\left(1-\frac{\kappa(\bx)(\xi^2+\eta^2)}{\rho(\bx)\omega^2}\right)^{-1/2}.
\end{equation}
Here $\xi$, $\eta$, and $\omega$ are the dual Fourier variables to
$x$, $y$, and $t$ respectively. The downgoing assumptions means that
for local planewave components of $P_sp$, the quantity inside the
square root is positive. Thus $\Lambda^{\pm}_{z_s}$ has real principal
symbol (in fact, the entire symbol is real) hence defines an
asymptotically symmetric operator:
\begin{equation}
  \label{eqn:lamappsim}
  (\Lambda^{\pm}_{z_s})^T \approx \Lambda^{\pm}_{z_s}.
\end{equation}
(For more on this, see \cite{StefUhl:05}.)
The analysis also reveals that the solution components not continuous
at $z=z_s$ are odd there:
\begin{equation}
  \label{eqn:odd1}
  \lim_{z\rightarrow z_s^+} v^{\pm}_{z} \approx - \lim_{z\rightarrow z_s^-}
  v^{\pm}_{z}
\end{equation}
for the solution of \ref{eqn:awepm} with $f_s=0$.
Similarly, 
\begin{equation}
  \label{eqn:odd2}
  \lim_{z\rightarrow z_s^+} p^{\pm}\approx - \lim_{z\rightarrow z_s^-}
  p^{\pm}
\end{equation}
for the solution of \ref{eqn:awepm} with $h_s=0$. Here ``$\approx$''
means in the sense of high frequency asymptotics, that is, that the
difference between the two sides is relatively smooth, hence small if
the data is highly oscillatory. Therefore if $f_s=0$ in system \ref{eqn:awepm},
\begin{equation}
  h_s = \Lambda^{\pm}_{z_s}P_sp^{\pm} = -[v^{\pm}_{z}]|_{z=z_s} \approx -2
  P_sv^{\pm}_{z}
  \label{eqn:tracejump10}
\end{equation}
Similarly, if $h_s=0$ in system \ref{eqn:awepm}, then
\begin{equation}
  \label{eqn:tracejump20}
  f_s = -[p^{\pm}]|_{z=z_s} \approx -2 P_s p^{\pm}.
\end{equation}
Thus $f_s$ determines approximately the boundary value of $p^{\pm}$,
as a solution of the acoustic wave system in the half-space
$z>z_s$. However, as repeated in equation \ref{eqn:tracejump10}, a
solution with this boundary value is also the restriction to $z>z_s$
of a solution to \ref{eqn:awepm} with $f_s=0$ and $h_s=
\Lambda^{\pm}_{z_s}P_sp^{\pm}$. Therefore if
\begin{equation}
  \label{eqn:hfcondn}
  h_s =-\frac{1}{2}\Lambda^{\pm}_{z_s}f_s,
\end{equation}
then the pressure boundary value $P_sp^{\pm}$ is the
same for the solutions of \ref{eqn:awepm} for source vectors $(h_s,0)$
and $(0,f_s)$. Since the pressure boundary values are the same, the solutions
in $z>z_s$ are the same. In particular, since $z_r>z_s$ and ${\cal
  S}^{\pm}_{z_s,z_r}(h_s,f_s)^T = (P_rp^{\pm},P_rv^{\pm}_z)^T$, it follows
that
\begin{equation}
  \label{eqn:snull}
  {\cal S}^{\pm}_{z_s,z_r}\left(\frac{1}{2}\Lambda^{\pm}_{z_s}f_s,f_s\right)^T \approx 0.
\end{equation}

Equation \ref{eqn:snull} states the relation between the columns of $
{\cal S}^{\pm}_{z_s,z_r}$ mentioned in the introduction to this
section.

\section{Time Reversal}

Recall that the source vector $(h_s,f_s)$ is assumed to produce a
downgoing field $(p^+,\bv^+)$, that is, emanates high-frequency energy only along
rays that make an angle with the vertical bounded below by a common
minimum angle. Such rays leave $\Omega$ within a common maximum
time. Consequently (Appendix B), in the
slab $z_s<z<z_r$, the field $(p^+,\bv^+)$ approximates the solution of an
anti-causal evolution equation. Choose $\chi(t)$ to be a smooth function
that is $= 0$ for $t \gg 0$ and $=1$ at times when near rays carrying
high-frequency energy in $(p^+,\bv^+)$ cross $z=z_r$. Define 
$(\tilde{p}^-,\tilde{\bv}^-)$ to be the solution in the half-space
$\Omega \times \bR$ of
\begin{eqnarray}
\label{eqn:revawe}
  \frac{1}{\kappa}\frac{\partial \tilde{p}^-}{\partial t} & = & - \nabla \cdot \tilde{\bv}^-, \nonumber \\
  \rho\frac{\partial \tilde{\bv}^-}{\partial t} & = & - \nabla \tilde{p}^-,\nonumber \\
  \tilde{p}^- & =& 0,  \mbox{ for } t \gg 0\\ 
  \tilde{\bv}^- & = & 0 \mbox{ for } t \gg 0\\
  P_r\tilde{p}^- &=& \chi P_rp^+ . 
\end{eqnarray}
That is, $\tilde{p}^-$ has the same boundary value on $z=z_r$ as
$p^+$, except for low-frequency residue that is muted by
$\chi$. Therefore
$p^+ \approx \tilde{p}^-, \bv^+ \approx \tilde{\bv}^-$ near
$z=z_r$. Since the right-hand sides in the system \ref{eqn:awepm} are
singular only on $z=z_s$, and the high-frequency components of
$(p^+,\bv^+)$ are carried by downgoing rays, these differ negligibly
from the the high-frequency components of
$(\tilde{p}^-,\tilde{\bv}^-)$ in the space-time slab $z_s<z<z_r$, and
the approximation holds throughout this region. In particular
$P_sv^+_z \approx P_s \tilde{v}^-_z$. In view of the relation
\ref{eqn:tracejump10},
\begin{equation}
  \label{eqn:tildevtohsubs}
  -2P_s\tilde{v}^-_z \approx h_s,
\end{equation}
so solution
of the system \ref{eqn:revawe} followed by restriction to $z=z_s$ and
multiplication by $-2$ 
approximately inverts the map $S^+_{z_s,z_r}: h_s \mapsto P_rp^+$.

Next observe that in view of the relation \ref{eqn:tracejump20}, and
the downgoing nature of the ray system carrying the high frequency
energy in $(p^+,\bv^+)$, the field $(\tilde{p}^-,\tilde{\bv}^-)$ is
actually the restriction to $z<z_r$ of the anti-causal solution of \ref{eqn:awepm}
with $z_s$ replaced by $z_r$, zero constitutive defect, and vertical
load given by the jump in pressure at $z=z_r$ - for this field, use
the same notation. Continuity of vertical
velocity $\tilde{v}^-_z$ at $z=z_r$ implies that the vertical load is
\[
  f_r = -[\tilde{p}^-]|_{z=z_r} =-(\lim_{z\rightarrow
    z_r^+}\tilde{p}^- - \lim_{z\rightarrow
    z_r^-}\tilde{p}^-)
\]
\[
  \approx 2 P_r \tilde{p}^- = 2 P_r p^+
\]
(from the definition \ref{eqn:defsamp}, $P_r$ is the limit from the
left). Thus
\[
  P_s \tilde{v^-_z} \approx V^-_{z_r,z_s}(2 P_rp^+) \approx
  2V^-_{z_r,z_s}S^+_{z_r,z_s}h_s.
\]
so
\[
  h_s \approx -2 P_s v^+_z \approx -2 P_s \tilde{v}^-_z \approx
  -4V^-_{z_r,z_s}S^+_{z_r,z_s}h_s
\]
Combine this observation with \ref{eqn:tildevtohsubs} to obtain
\[
 -4  V^-_{z_r,z_s} S^+_{z_s,z_r}  \approx  I,
\]
This relation combines with the identity \ref{eqn:trtrcomp} to
yield the first main result of this section:
\begin{eqnarray}
  \label{eqn:approxinv}
  (V^+_{z_s,z_r})^T S^+_{z_s,z_r} & \approx & \frac{1}{4}I, \nonumber\\
  (S^+_{z_s,z_r})^T V^+_{z_s,z_r} & \approx & \frac{1}{4}I, \nonumber\\
  V^+_{z_s,z_r} (S^+_{z_s,z_r})^T & \approx & \frac{1}{4}I, \nonumber\\
  S^+_{z_s,z_r} (V^+_{z_s,z_r})^T & \approx & \frac{1}{4}I.
\end{eqnarray}.
The second equation is simply the transpose of the first, and the
last two follow by by an exactly analogous argument using time
reversal and interchange of the roles of $z_s$ and$z_r$.

The conclusion is significant enough to merit restating in English:
provided that high-frequency energy in the various fields is carried
along downgoing ray fields, the transpose of $V^+$ is an approximate
inverse to $S^+$, modulo a factor of 4. To recover the pressure source
$h_s$ generating a pressure gather $P_rp$ at $z=z_r$, multiply the
latter by -2, then apply the transpose of $V^+_{z_s,z_r}$ to this
gather, reading out a vertical velocity field at $z=z_s$. Multiply
again by -2 and you have a high-frequency approximation to $h_s$.

\section{Unitarity}

The next chapter in this story recognizes the relations in display
\ref{eqn:approxinv} as asserting the approximate unitarity of
$S^+_{z_s,z_r}$.

The matrix identity \ref{eqn:snull} implies a relation between $S, V,$
and $\Lambda$ of some interest in itself. After minor re-arrangement, the second row of reads
\begin{equation}
  \label{eqn:snull2}
-\frac{1}{2}\Pi_1{\cal S}^{\pm}_{z_s,z_r}\Pi_0^T\Lambda^{\pm}_{z_s}  \approx
V^{\pm}_{z_s,z_r}.
\end{equation}
In these relations, the projection on the left picks out the vertical velocity component
of a downgoing wavefield at $z=z_r$: that is,
\[
-\frac{1}{2}\Pi_1{\cal S}^{\pm}_{z_s,z_r}\Pi_0^T\Lambda^{\pm}_{z_s}P_sp^+
=-\frac{1}{2}P_r v_z^+,
\]
where $(p^+, \bv^+)$ solve the system \ref{eqn:awepm} with $f_s=0$ and
$h_s = \Lambda^{\pm}_{z_s}P_sp^+$. On the other hand, from relation
\ref{eqn:tracejump10},
\[
  P_r v_z^+ = -\frac{1}{2}\Lambda^+_{z_r}P_r p^+
\]
where
\[
  P_r p^+ = \Pi_0{\cal S}^{+}_{z_s,z_r}\Pi_0^T\Lambda^{+}_{z_s}P_s
  p^+
\]
\[
  = S^+_{z_s,z_r}\Lambda^{+}_{z_s}P_sp^+
\]
Therefore combining the last two equations with \ref{eqn:snull2},
obtain
\begin{equation}
  \label{eqn:sv}
  \frac{1}{4}\Lambda^+_{z_r}S^+_{z_s,z_r}\Lambda^{+}_{z_s} = V^+_{z_s,z_r}.
\end{equation}
This is the promised relation.

As shown in the last section, $4(V_{z_s,z_r}^+)^T$ is approximately
inverse to $S^{+}_{z_s,z_r}$. Therefore, transposing both sides of
equation \ref{eqn:sv} and using \ref{eqn:approxinv}, obtain
\begin{equation}
  \label{eqn:almostunitary}
  4(V_{z_s,z_r}^+)^TS^+_{z_s,z_r} = [ (\Lambda^+_{z_s})^T
  (S^{+}_{z_s,z_r})^T(\Lambda^+_{z_r})^T]S^{+}_{z_s,z_r} \approx I.
\end{equation}

The remarkable feature of the identity \ref{eqn:almostunitary} is that
it exhibits an approximate right inverse of $S^+$ as an adjoint with
respect to a weighted inner product - or it would, if the operators
$(\Lambda^+)$ were symmetric positive definite. As noted earlier,
these operators are only approximately symmetric, though they are
positive semi-definite. That is not a great obstacle, however:
symmetrizing them in the obvious way commits a negligible error, of
the sort that this paper already neglects wholesale. That is,
\begin{equation}
  \label{eqn:unitary}
  [ \frac{1}{2}((\Lambda^+_{z_s})^T+ \Lambda^+_{z_s})
  (S^{+}_{z_s,z_r})^T \frac{1}{2}((\Lambda^+_{z_r})^T+
  \Lambda^+_{z_r})]S^{+}_{z_s,z_r} \approx I.
\end{equation}

The symmetrized $\Lambda$ operators are at least positive
semi-definite, hence define (at least) semi-norms.
Similar relations have been derived for other scattering operators,
and have been used to accelerate iterative solutions of inverse
scatering problems: \cite{DafniSymes:SEG18b} review some of this
literature.

\section{Accelerated Iterative Inversion}

For convenience, in this section write $S$ in place of
$S^+_{z_s,z_r}$. Also abbreviate the symmetrized $\Lambda$ operators
using notation suggesting weight
operators in model and data spaces:
\begin{eqnarray}
  W_m^{-1}&=& \frac{1}{2}((\Lambda^+_{z_s})^T+
              \Lambda^+_{z_s}),\nonumber \\
  W_d &=& \frac{1}{2}((\Lambda^+_{z_r})^T+ \Lambda^+_{z_r}).
          \label{eqn:wdef}
\end{eqnarray}
The identification of the symmetrized $\Lambda^+_{z_s}$ as the inverse
of another operator $W_m$ is formal, since the former operator is
likely to have null (or nearly-null) vectors due to aperture-related
amplitude loss. Since some version of $W_m$ is essential in the
formulation for effective preconditioning, I will derive a usable
candidate to stand in for it below.

Adopting Hilbert norms defined by the operators $W_m$ and $W_d$ in its
domain and range respectively, the adjoint of $S$ is given by
\begin{equation}
\label{eqn:wadj}
S^{\dagger} = W_m^{-1}S^TW_d,
\end{equation}

In this notation, the relation \ref{eqn:unitary} takes the form
\begin{equation}
  \label{eqn:wunitary}
  S^{\dagger}S \approx I.
\end{equation}
That is to say, $S$ is approximately unitary with respect to the
weighted norms defined by $W_m$ and $W_d$. Therefore a Krylov space
method employing these norms will converge rapidly, at least for the
well-determined components of the solution.

The most convenient arrangement the Conjugate Gradient (CG) algorithm
taking advantage of the structure \ref{eqn:wadj} is the {\em
  Preconditioned CG}. Allowing that the fit error will be measured by
the data space norm, the least squares problem to be solved is not
just $Sh \approx d$, but a regularized version:
\begin{equation}
  \label{eqn:einv}
  \mbox{minimize}_h \|Sh-d\|^2_d + \alpha^2 \|Ah\|^2_m
\end{equation}

\noindent {\bf Remark:} recall that the modified data space norm $\|d\|_d^2 = \langle
d, W_d d\rangle$ has physical meaning: for acoustics, it is
proportional to the power transmitted to the fluid by the source.

The minimizer of the objective defined in equation \ref{eqn:einv}
solves the normal equation
\begin{equation}
  \label{eqn:norm0}
  (S^{\dagger}S + \alpha^2 A^{\dagger}A)h = S^{\dagger}d 
\end{equation}
where the weighted adjoint $S^{\dagger}$ has already been defined in equation \ref{eqn:wadj}, and $A^{\dagger}$ is the adjoint of $A$ in the weighted model space norm defined by $W_m$, namely
\begin{equation}
  \label{eqn:aadj}
  A^{\dagger} = W_m^{-1}A^TW_m.
\end{equation}

Note that the normal operator appearing on the left-hand side of
\ref{eqn:norm0} is not an approximate identity, due to the presence of
the regularization term: the spectrum increases in spread with
increasing $\alpha$, leading to slower convergence. Fortunately for
the present setting, the operators $W_m^{-1}$, $A$, and $W_m$
approximately commute (they are scalar {\em pseudodifferential}, once
the difficulties with the definition of $W_m$, mentioned above, are
taken care of). Scalar pseudodifferential operators approximately
commute, so $A^{\dagger} \approx A^T$. Therefore
\begin{equation}
  \label{eqn:normapprox}
  S^{\dagger}S + \alpha^2 A^{\dagger}A \approx I + \alpha^2A^TA
\end{equation}
Recall that $A$ is simply multiplication by the Euclidean distance to
the physical source point $\bx_s$: $A u (\bx) = |\bx-\bx_s|u(\bx),
A^TAu(\bx) = |\bx-\bx_s|^2u(\bx)$. So the equation $(I+\alpha^2
A^TA)u=b$ is trivial to solve, and this is a key characteristic of a
good preconditioner. However this observation must be combined with
the weighted norm structure.

Rewrite the normal equation \ref{eqn:norm0} as
\begin{equation}
  \label{eqn:norm1}
  W_m^{-1}(S^TW_dS + \alpha^2 A^TW_mA)h = W_m^{-1}S^TW_md 
\end{equation}
Since $W_m$ is self-adjoint and positive semidefinite, the common factor on both sides of \ref{eqn:norm1} can be re-written as
\begin{equation}
  \label{eqn:normpart}
  Nh = (S^*S + \alpha^2 A^*A)h = S^*d 
\end{equation}
in which $S^*, A^*$ are the adjoints with the original (Euclidean)
inner product in the domains but the weighted inner product in data
space:
\begin{eqnarray}
  \label{eqn:sadjwt}
  S^* &=& S^T W_d,\\
  A^* &=& A^T W_m.
\end{eqnarray}
Note the $S^*S$ and $A^*A$ are symmetric in the Euclidean sense, so
equation \ref{eqn:normpart} is a symmetric positive (semi-)definite
linear system, just the sort of thing for which the 
The Preconditioned Conjugate Gradient (``PCG'') algorithm was
designed. PCG for solution
of equation \ref{eqn:normpart} with preconditioner $M$ is usually
written as Algorithm 1 (see for example \cite{Golub:2012}):

\begin{algorithm}[H]
\caption{Preconditioned Conjugate Gradient Algorithm, Standard Version}
\begin{algorithmic}[1]
\State Choose $h_0=0$ 
  \State $r_0 \gets S^*d$
  \State $p_0 \gets M^{-1} r_0$
  \State $g_0 \gets p_0$
  \State $q_0 \gets Np_0$
  \State $k \gets 0$
  \Repeat
  \State $\alpha_k \gets \frac{\langle g_k,r_k \rangle}{\langle p_k,q_k\rangle}$
  \State $h_{k+1} \gets h_k + \alpha_k p_k$
  \State $r_{k+1} \gets r_k - \alpha_kq_k$
  \State $g_{k+1} \gets M^{-1} r_{k+1}$
  \State $\beta_{k+1} \gets \frac{\langle g_{k+1},r_{k+1}\rangle}{\langle g_k,r_k\rangle}$
  \State $p_{k+1}\gets g_{k+1}+\beta_{k+1}p_k$
  \State $q_{k+1} \gets Np_{k+1}$
  \State $k \gets k+1$
  \Until{Error is sufficiently small, or max iteration count exceeded} 
\end{algorithmic}
\end{algorithm}
The iteration converges rapidly if $M^{-1}N \approx I$. This is true
if and only if the symmetrized operator $M^{-1/2}NM^{-1/2} \approx I$,
which is in turn true if the eigenvalues of $M^{-1/2}NM^{-1/2}$ are
close to 1 (actually works well is most of these eigenvalues are close
to 1, and the rest are small - which is the case for the current
problem)..  Further, PCG is computationally effective is M is easy to
invert.

From \ref{eqn:normapprox} and \ref{eqn:norm1}, it follows that
\[
  W_m^{-1}(S^TW_dS + \alpha^2 A^TW_mA) \approx I + \alpha^2 A^TA.
\]
This observation suggests using $M=W_m(I+\alpha^2A^TA)$. This choice
is not symmetric, but since the operators on the right-hand side are
scalar pseudodifferential hence commute, it is equivalent to use of
\begin{eqnarray}
  M         &=&(I+\alpha^2A^TA)^{1/2}W_m(I+\alpha^2A^TA)^{1/2},\nonumber \\
  M^{-1}
            &=&(I+\alpha^2A^TA)^{-1/2}W_m^{-1}(I+\alpha^2A^TA)^{-1/2}.
                \label{eqn:defprecond}
\end{eqnarray}
With this choice, \ref{eqn:normapprox} implies that
$M^{-1}N \approx I$, also $M$ is symmetric. As already mentioned,
powers of $I + \alpha^2A^TA$ are trivial to compute, given the choice
of $A$ made here. We will examine fast algorithms for computing
$W_m^{-1}$ = the symmetrized pressure-to-source operator in the next
section. Note that only $M^{-1}$, hence only $W_m^{-1}$, appears in
Algorithm 1.

\section{Computing and Symmetrizing $\Lambda$}

Computations of $\Lambda^{\pm}_{z_s}$ and its transpose are clearly critical steps in an
implementation of the PCG algorithm outlined in the preceding section.
Direct computation of the pressure-to-source operator $\Lambda^{\pm}_{z_s}$, for instance by solving
\ref{eqn:awepm} and reading off $P_sv^{\pm}_{z}$, turns out to be
numerically ill-behaved. The relation \ref{eqn:snull} provides and
alternative approach, taking advantage of the accurate approximate
inverse to $S^+_{z_s,z_r}$ constructed above. The first row of
\ref{eqn:snull}, slightly rearranged, is
\begin{equation}
  \label{eqn:lamidea}
  \Pi_0{\cal S}^+_{z_s,z_r}\Pi_1^T f_s \approx  -\frac{1}{2} S^+_{z_s,z_r}\Lambda^+_{z_s}f_s.
\end{equation}
The approximate inverse construction for $S^++_{z_s,z_r}$ permits (approximate)
solution of this equation for $\Lambda^+_{z_s}f_s$: apply
$4(V^{+}_{z_s,z_r})^T$ to both sides of equation \ref{eqn:lamidea} and
use the first equation in the list \ref{eqn:approxinv} to get
\begin{equation}
  \label{eqn:lamident}
  \Lambda^+_{z_s} \approx -8(V^{+}_{z_s,z_r})^T\Pi_0 {\cal S}^+_{z_s,z_r}\Pi_1^T.
\end{equation}
This identity is the major result of this section: it shows how
to compute that action of $\Lambda^+_{z_s}$ by propagating the input
pressure trace, identified as a source for the velocity evolution,
forward in time from $z_s$ to $z_r$
reading off the pressure trace on $z=z_r$, identifying it once more as
a point load (source for velocity), propagating it backwards in time from $z_r$ to
$z_s$, and finally reading off the velocity trace, interpreted as a
pressure evolution source on $z_s$. 

The importance of this result lies in the failure of the obvious
method for computing the action of $\Lambda^{\pm}_{z_s}$, namely to
employ the pressure trace as a source in the velocity equation ($f_s$,
in the notation used above) at $z=z_s$, and read off the velocity
field also at $z=z_s$. This difficulty is related to the existence of
tangentially propagating waves and the lack of continuity of the trace
operator. The method implicit in equation \ref{eqn:lamident} avoids
this difficulty by propagating the fields a positive distance in $z$:
assuming as always that the causal fields are downgoing, this step
eliminates any tangentially propagating fields from consideration.

A deeper study of the pressure-to-source operator (or of the closely
related Dirichlet-to-Neumann operator for the second order wave
equation, see \cite{StefUhl:05}) shows that it is approximately
dependent only on the model coefficients near the source surface
($z=z_s$ in this case). Since the homogeneous and lens models are
identical near this surface, it is unsurprising that these figures are
very close to the previous two.  However an even more useful
observation is that the calculations in the approximation
\ref{eqn:lamident} could just as well be carried out in a much smaller
region around the source surface, and produce a result that is
functionally identical in that it will serve as a source for the same
acoustic fields globally, with small error. In effect, equation
\ref{eqn:lamident} involving propagation from source ($z=z_s$) to
receiver ($z=z_r$) surfaces is altered by replacing $z_r$ with a
receiver datum $z_s+\Delta z$ considerably closers to $z_s$:
\begin{equation}
  \label{eqn:lamnear}
  \Lambda^+_{z_s} \approx -8(V^{+}_{z_s,z_s+\Delta z})^T\Pi_0 {\cal
    S}^+_{z_s,z_s+\Delta z}\Pi_1^T.
\end{equation}
Using a receiver datum closer to the source surface has two favorable consequences:
\begin{itemize}
\item The computational domain can be smaller than is necessary to
  simulate the target data, as it need only contain the source
  surface and the receiver datum implicit in
  equation \ref{eqn:lamident}. This shrinkage of the computational
  domain can lead to substantial improvements in computational
  efficiency.
\item Since the receiver data may be chosen much closer to the
  source surface that is the case for the target data, the effective
  aperture active in the relation \ref{eqn:lamident} can be much
  larger, producing an estimated source gather much less affected by
  aperture limitation.
\end{itemize}

As mentioned in the last section, computation of the transpose of
$\Lambda^+$ (exact, not approximate in the high frequency sense) is
critical to the successful construction of the preconditioner. The
relation \ref{eqn:lamident} does not provide a computation for this
operator. However set
\begin{equation}
  \label{eqn:lamtilde}
  \tilde{\Lambda}^+_{z_s} = -8(V^{+}_{z_s,z_s+\Delta z})^T\Pi_0 {\cal
    S}^+_{z_s,z_s+\Delta z}\Pi_1^T.
\end{equation}
Then \ref{eqn:lamnear} can be rewritten
\[
  \Lambda^+_{z_s} \approx \tilde{\Lambda}^+_{z_s}.
\]
Of course, all of the examples so far show images of
$\tilde{\Lambda}^+_{z_s}$.

Since successful preconditioning requires only approximate inversion,
use of $\tilde{\Lambda}^+_{z_s} $ in place of $\Lambda^+_{z_s}$ will
still yield a working preconditioner, and the former can be transposed
to machine precision via the definition \ref{eqn:lamtilde} and the adjoint state
method (equations \ref{eqn:trtr} \ref{eqn:trtrcomp}):
\begin{equation}
  \label{eqn:lamtransp}
  (\tilde{\Lambda}^+_{z_s})^T = -8 \Pi_1 ({\cal S}^+_{z_s,z_s+\Delta
    z})^T\Pi_0^T V^{+}_{z_s,z_s+\Delta z}
\end{equation}

The model space weight operator $W_m^{-1}$ introduced in the last section is
replaced by its asymptotic approximation
\[
\frac{1}{2}(\tilde{\Lambda}^+_{z_s} +
  (\tilde{\Lambda}^+_{z_s})^T)
\]
\[
  \approx -8 \left((V^{+}_{z_s,z_s+\Delta z})^T\Pi_0 {\cal
    S}^+_{z_s,z_s+\Delta z}\Pi_1^T \right.+
\]
\[
  \left. \Pi_1 ({\cal S}^+_{z_s,z_s+\Delta
      z})^T\Pi_0^T V^{+}_{z_s,z_s+\Delta z}\right)
\]
\begin{equation}
  \label{eqn:winvcomp}
  = -4(\Pi_1 ({\cal S}^+_{z_s,z_s+\Delta z})^T(\Pi_0^T\Pi_1 +
  \Pi_1^T\Pi_0) {\cal S}^+_{z_s,z_s+\Delta z}\Pi_1^T) = \tilde{W}_m^{-1} .
\end{equation}
with a similar definition for the replacement $\tilde{W}_d$ of $W_d$.

This identity shows that only one forward and one adjoint simulation
are necessary to compute the action of $\tilde{W}_{m,d}$. The operator
in the center of the expression on the right-hand side, $\Pi_0^T\Pi_1
+ \Pi_1^T\Pi_0$, simply exchanges the components of the acoustic
fields, passing the velocity field as a pressure source and the
pressure field as a velocity source.

One more computation is required for the full implementation of the
preconditioning strategy explained in the last section: $W_m$ is
required, not just $W_m^{-1}$. Note that $W_m$ plays two roles in the
second term in equation \ref{eqn:norm1}: it is the weight matrix for
both the domain and range norms for $A$. It is perfectly OK for one of
these to be replaced by an asymptotic approximation, so long as it is
symmetric and computable (and at least semi-definite). The second row
in equation \ref{eqn:snull} appears as \ref{eqn:snull2} above:
introducing (formally) the inverse of $\Lambda^+$,
\begin{equation}
  \label{eqn:snull2mod}
-\frac{1}{2}\Pi_1{\cal S}^{+}_{z_s,z_r}\Pi_0^T  \approx
V^{+}_{z_s,z_r}(\Lambda^{+}_{z_s})^{-1}
\end{equation}
whence from the second line in display \ref{eqn:approxinv}
\begin{equation}
  \label{eqn:laminv}
-\frac{1}{8}(S^+_{z_s,z_r})^T\Pi_1{\cal S}^{+}_{z_s,z_r}\Pi_0^T  
\approx (\Lambda^{+}_{z_s})^{-1}
\end{equation} 
and
\begin{equation}
  \label{eqn:laminvtr}
-\frac{1}{8}\Pi_0({\cal S}^{+}_{z_s,z_r})^T\Pi_1^T S^+_{z_s,z_r}
\approx ((\Lambda^{+}_{z_s})^{-1})^T.
\end{equation}
Using the definition \ref{eqn:sdef} of $S^+_{z_s,z_r}$, the
symmetrized $\Lambda^{-1}$ is
\begin{equation}
  \label{wcomp}
 \tilde{W}_m = -\frac{1}{16}\left(\Pi_0 ({\cal S}^{+}_{z_s,z_r})^T
   (\Pi_1^T \Pi_0 + \Pi_0^T \Pi_1) {\cal S}^{+}_{z_s,z_r} \Pi_0^T\right)
   \approx \frac{1}{2}((\Lambda^{+}_{z_s})^{-1}
   +((\Lambda^{+}_{z_s})^{-1})^T) .
\end{equation}
Comparison with the definition \ref{eqn:winvcomp} shows that
$\tilde{W}_m$ and $\tilde{W}_m^{-1}$ differ only in the initial and
final projection factors (and overall scale), and in particular either
can be computed for the cost of a forward/adjoint operator pair. Note
that $\tilde{W}_m^{-1}$ is inverse to $\tilde{W}_m$ only in an
approximate (asymptotic, aperture-limited) sense.

\section{Numerical Examples}
This section illustrates the most important conclusions developed in
the preceding sections by finite difference wavefield simulation.

\subsection{Synthetic models and simulation}
To illustrate the structure described in the preceding section, I
introduce two 2D acoustic models, one spatially homogeneous, the other
highly refractive. The first, homogenous model has $\kappa = 4$
GPa and $\rho = 1$ g/cm$^3$ throughout a rectangular domain of size 8 km ($x$) $\times$ 4 km
($z$). The second, refractive, model is a perturbation of the first by
a low-velocity acoustic lens positioned in the center of the
rectangle (Figure \ref{fig:bml0}. To produce this structure, the density is chosen
homogeneous as in the first model, while the bulk modulus decreases to
from 4 GPa outside the lens to 1.6 GPA in its center, as shown in Figure \ref{fig:bml0}.

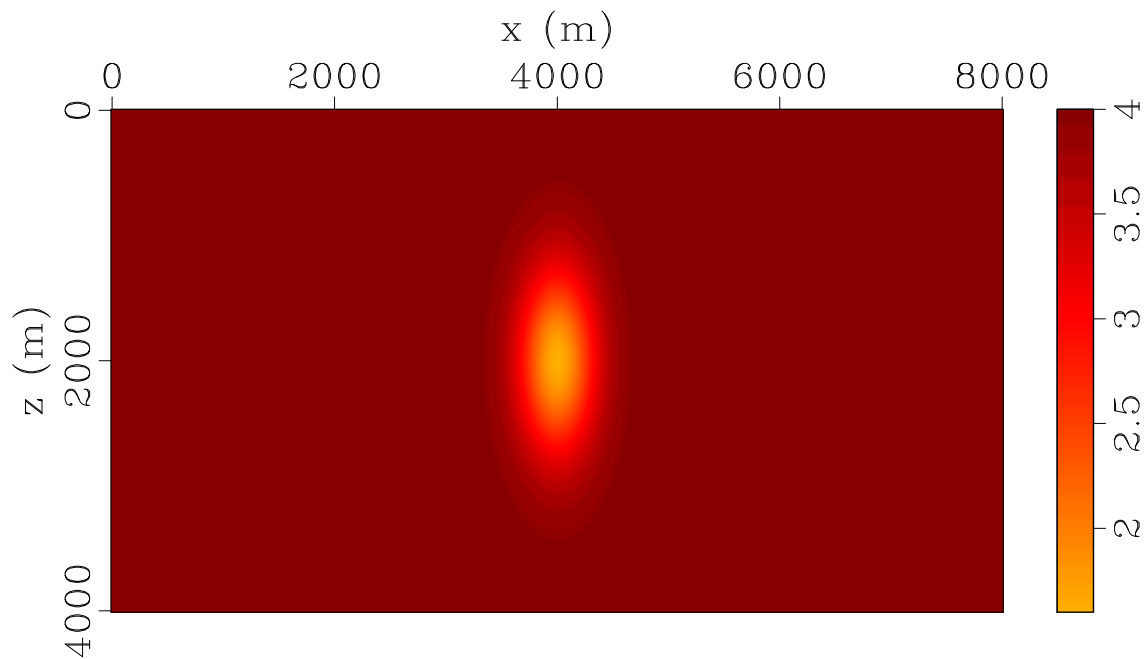
\plot{bml0}{width=\textwidth}{Bulk modulus, lens model. Color scale is 
in GPa. Positions of point source and receiver line indicated.}

Discretization is conventional, with a rectangular grid and staggered
finite difference scheme \cite[]{Vir:84} of order 2 in time and
2$k$ in space; for most of the experiments reported below,
$k=4$. Absorbing boudary conditions of perfectly matched layer type
are applied at all boundaries of the simulation rectangle \cite[]{Habashy:07}.
Sampling operators such as $P_r$ are implemented via linear
interpolation, and source insertion via adjoint linear interpolation
(as noted above, in the continuum limit, sources are represented via
adjoint sampling). Steps in $x$ and $z$ are the same. In the following
examples, $\Delta x = 20$ m. This choice limits the temporal frequency
of accurately computed fields to rougly 12 Hz.

\cite{GeoPros:11} gives a description of the code
implementation, out-of-date in a few respects but overall accurate.
The implementation uses the discrete adjoint state method and
auto-generated code \cite[]{TapenadeRef13}, to assure that the
computed adjoint operators are adjoint at the level of machine
precision to the computed operators. The reverse-time storage issue is
resolved through the optimal checkpointing technique
\cite[]{Griewank:book,Symes:06a-pub}, again without loss of
precision. This procedure results in computed adjoints for
$S^+_{z_s,z_r}$ and other operators that pass usual test for adjoint
accuracy, comparing inner products with pseudorandom input vectors,
with errors well under machine precision.



The horizontal line of receivers sits at depth $z_r = $ 1000 m, the
(extended) sources at $z_s=3000$ m.  
Source and receiver $x$ ranges from $2000$ to
$6000$ m. Note that we have
reversed the order relation between $z_s$ and $z_r$ described in the
text ($z_s<z_r$). This difference is immaterial for the purpose of
illustrating the mathematical structures developed in the preceding paragraphs.




\subsection{Creating downgoing fields}
The downgoing condition constrains high-frequency energy of localized
plane wave components, hence could be enforced by dip>
filtering. However, a simpler approach is to construct fields that
must be entirely downgoing at the source and receiver surfaces by
virtue of ray geometry.

Note that a point source
on $z=z_s$ creates high frequency energy traveling on rays parallel
and nearly parallel to $z=z_s$, so that won't do. However, placing a
point source at a depth $z_d<z_s$ will work. Since the examples used here are
homogeneous in $z<z_s$, and the sampling region for extended sources
is a finite interval, all rays carrying high frequency energy cross
the source surface $z=z_s$ at a postive angle, and the field and its
traces are {\em a priori} downgoing. The same is obviously true at the
receiver surface for the homogeneous model, but is also true for the
lens model, as no rays are refracted horizontally {\em at the receiver surface}.

The choice of a point source at
$z_d=3500$ m, $x_d=3500$ m,  bandpass filter wavelet with
corner frequencies $1, 2.5, 7.5, 12.5$ Hz, gives the causal pressure and
velocity gathers at $z=z_s=3000$ m
shown in Figures \ref{fig:dsrcphh0} and \ref{fig:dsrcvzhh0}.  Since the
mechanical parameters in the homogeneous and lens models are the same
for $z<z_s$, and no rays return to this zone in either model, these
data are asymptotically the same for both models, and I show only the
homogenous medium results.

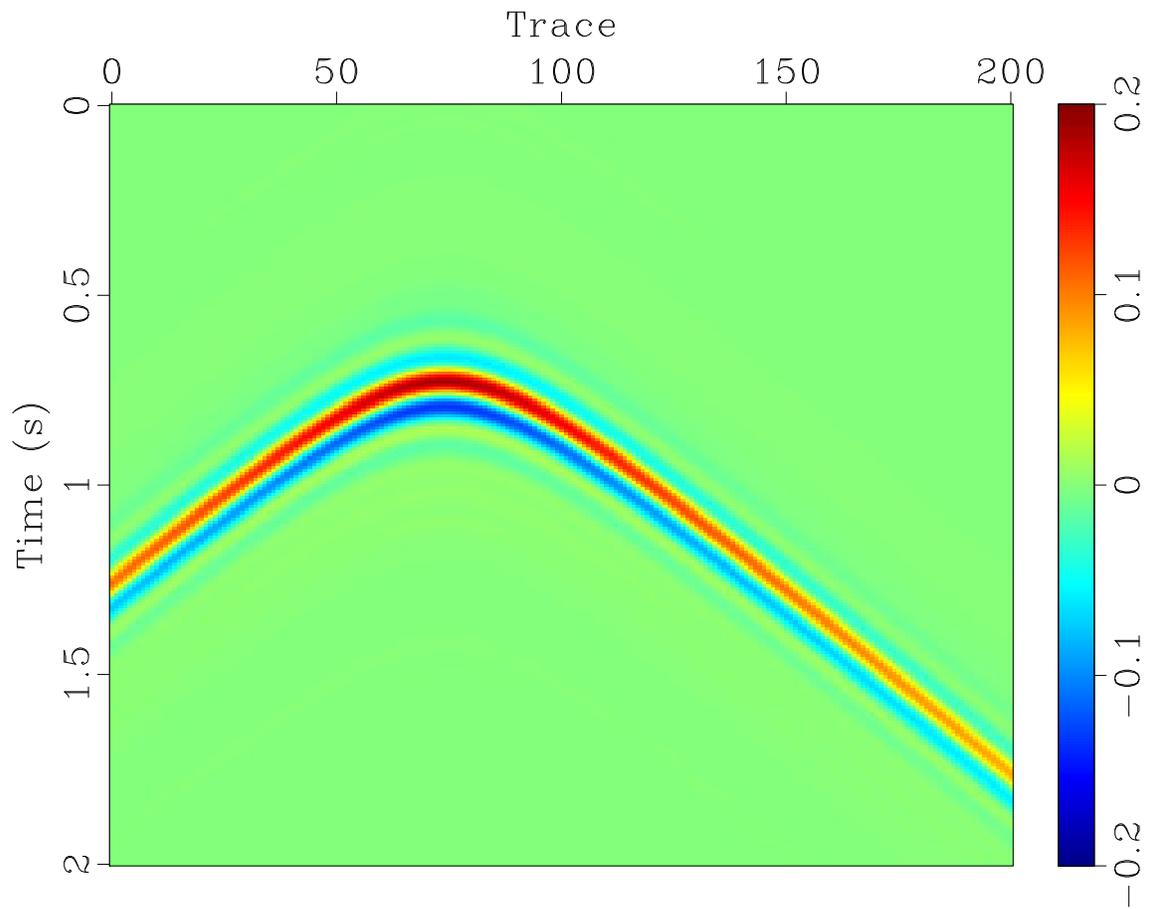
\plot{dsrcphh0}{width=\textwidth}{Trace $P_sp^+$ on $z=z_s=3000$ m of
  pressure field from point source at $z_d=3500$ m, $x_d=3500$ m,
  bandpass filter source. }

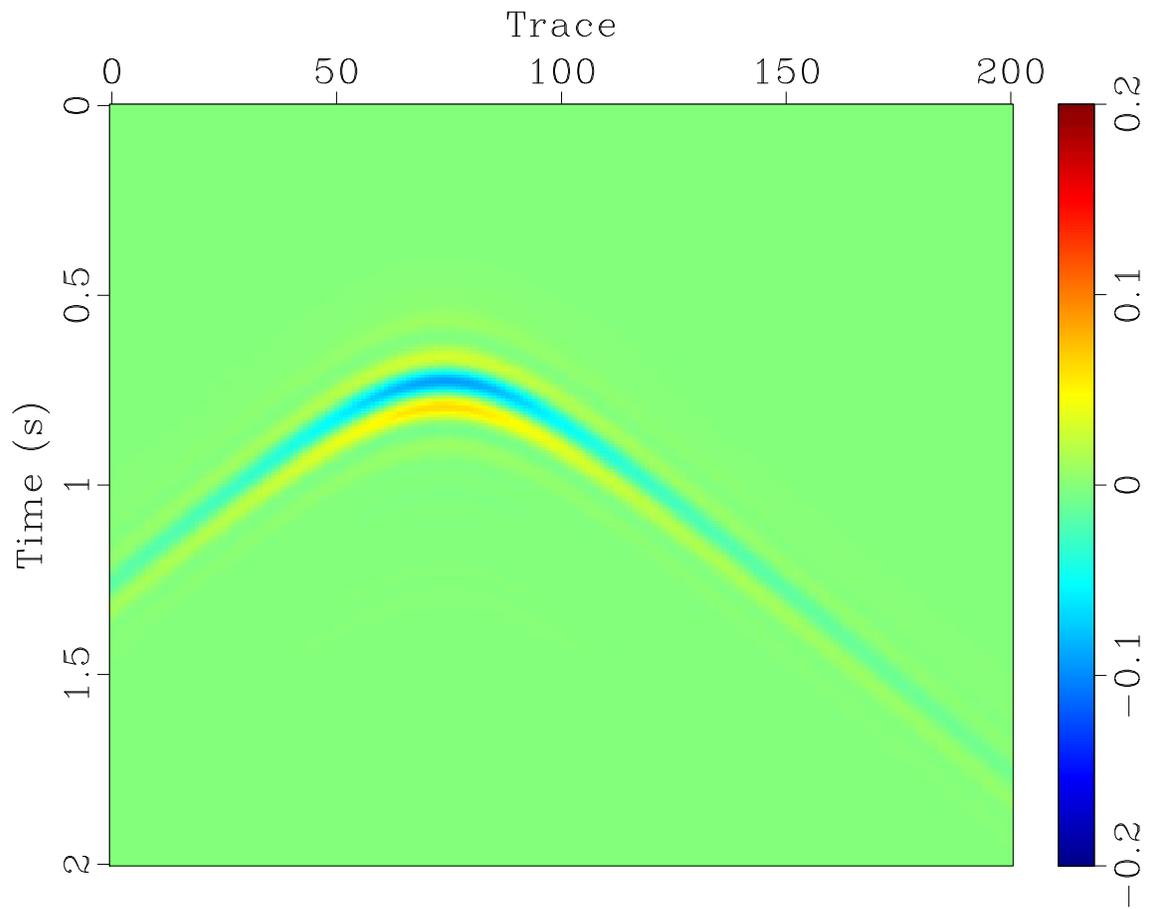
\plot{dsrcvzhh0}{width=\textwidth}{Trace $P_sv_z^+$ on $z=z_s=3000$ m of
  vertical velocity field from point source at $z_d=3500$ m, $x_d=3500$ m,
  bandpass filter source.}

\subsection{Equivalence of pressure, velocity sources}
These gathers are the pressure and
velocity traces
$(P_sp^+,P_sv^+_z)$ on $z=z_s$ of a downgoing acoustic field in
$z<z_s$, hence related by the operator $\Lambda^+_{z_s}$.
Equations \ref{eqn:snull}, \ref{eqn:tracejump10} and
\ref{eqn:tracejump20} show that these differ by a factor of -2 from
source functions $f_s$ and $h_s$ in the system \ref{eqn:awepm},
with $h_s=0$ and $f_s=0$ respectively, that generate the same acoustic
field in $z<z_s$, and in particular the same receiver traces on
$z=z_r$, at least asymptotically.


Figures
\ref{fig:drecplh0}, \ref{fig:dfwdplh0}, \ref{fig:daltplh0} show the pressure
gathers extracted at $z_r=1000$ m for the point source at $z=z_d$ and
for the two choices of extended source at $z=z_s$, on the same color
scale. The obvious similarity between the fields generated by the two
extended sources,
predicted by equation \ref{eqn:snull}, is confirmed by trace
comparisons in figures 
\ref{fig:drecplh0tr81},\ref{fig:daltplh0tr81}. Other traces are
equally similar.

\plot{drecplh0}{width=\textwidth}{Pressure gather at receiver depth 
  $z_r=1000$ m from field generated by causal solution of acoustic 
  system \ref{eqn:awepm} in the lens model described in the 
  text, with point pressure source (constitutive defect) at $z_d=3500$
  m, $x_d=3500$ m.}

\plot{dfwdplh0}{width=\textwidth}{Pressure gather at receiver depth 
  $z_r=1000$ m from field generated by causal solution of acoustic 
  system \ref{eqn:awepm} in the lens model described in the 
  text, with extended pressure source (constitutive defect) on
  $z=z_s=3000$ m given bythe field depicted in Figure
  \ref{fig:dsrcvzhh0} scaled by -2 
   ($h_s=-2P_sv_z^+=\Lambda^+_{z_s}P_sp^+$) and zero velocity source (vertical load) 
  ($f_s=0$).}

\plot{daltplh0}{width=\textwidth}{Pressure gather at receiver depth 
  $z_r=1000$ m from field generated by causal solution of acoustic 
  system \ref{eqn:awepm} in the lens model described in the 
  text, with extended velocity source (vertical load) on $z=z_s=3000$
  m given by the field depicted in Figure \ref{fig:dsrcphh0} scaled by -2
  ($f_s=-2P_sp^+$) and zero pressure source (constitutive defect)
  ($h_s=0$).}

\plot{drecplh0tr81}{width=\textwidth}{Overplot of traces 81 ($x=3600$) 
  from gathers shown in \ref{fig:drecplh0} (blue), \ref{fig:dfwdplh0}
  (red).}

\plot{daltplh0tr81}{width=\textwidth}{Overplot of traces 81 ($x=3600$)
  from gathers shown in \ref{fig:drecplh0} (blue), \ref{fig:daltplh0}
  (red).}

\subsection{Inversion by time reversal}
I have applied the approximate inversion procedure suggested in
equation \ref{eqn:approxinv} to the pressure gather shown in Figure
\ref{fig:drecplh0}, generated by a point source
at $z_d=3500$ m, $x_s=3500$ m, propagating in the 
lens
model (Figure \ref{fig:bml0}). I choose this example for two
reasons. First, the success of the inversion demostrates the
insensitivity of the time reversal method to ray multipathing
(triplication), evident in the data (Figure \ref{fig:drecplh0}).
Second, I will invert this data in the homogeneous model, that is,
construct sources that (approximately) reproduce the data using a
different material model than the one in which it was produced. This
capability is critically important in the application of the
approximate inversion in nonlinear inversion, where the early iterations
involve solution of the source estimation problem \ref{eqn:esis} at
(possibly very) wrong material models ${\bf c}$. Successful extension
methods maintain data fit throughout the course of the inversion.

As noted earlier, the acoustic field in this example is downgoing
throughout the simulation range. It can be
regarded as the result of either pressure or velocity source at $z=z_s$: the
pressure source gather $h_s$ (Figure \ref{fig:dhshh0})  is -2 times the
vertical velocity gather depicted in Figure \ref{fig:dsrcvzhh0}, the
velocity source gather $f_s$ (Figure \ref{fig:dfshh0}) is -2 times the
pressure gather depicted in Figure \ref{fig:dsrcphh0}.

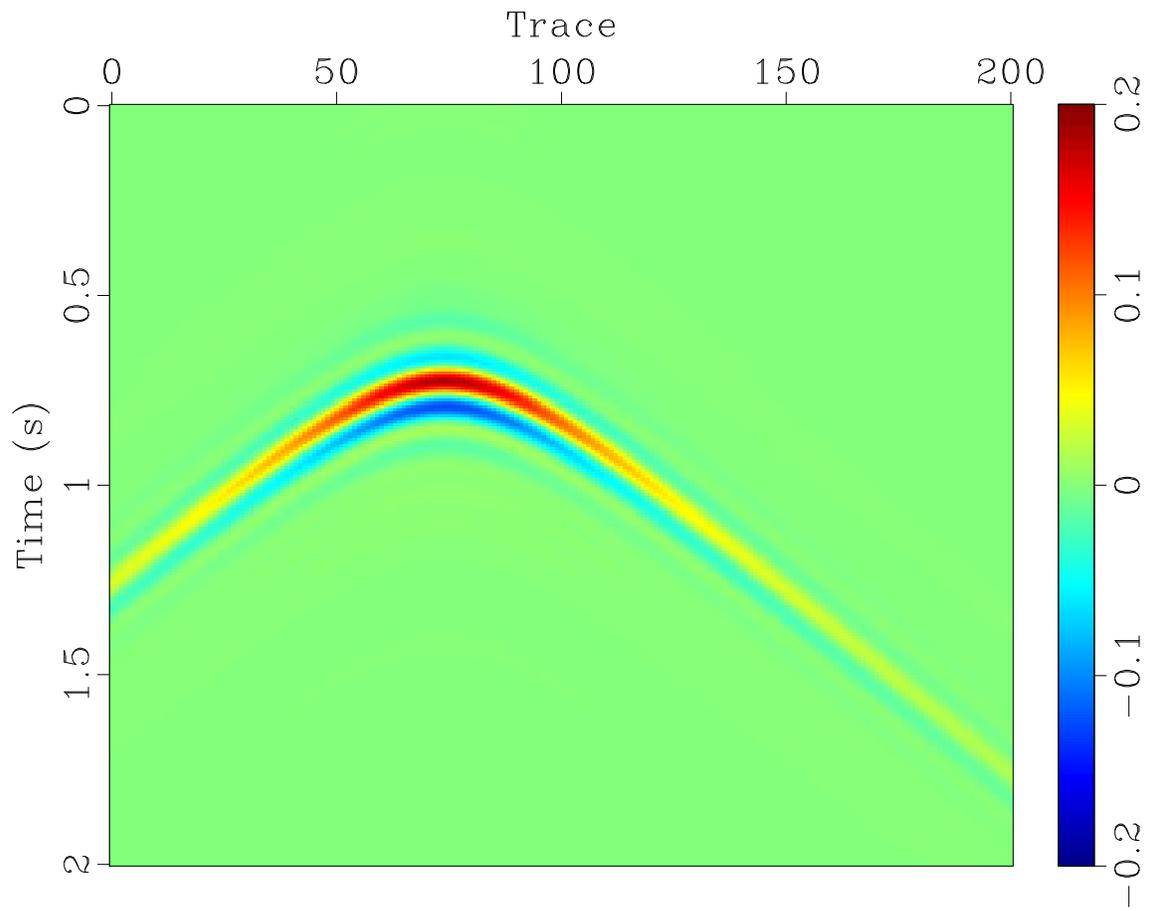
\plot{dhshh0}{width=\textwidth}{Pressure source gather = -2 $\times$
  vertical velocity gather (Figure \ref{fig:dsrcvzhh0}).}

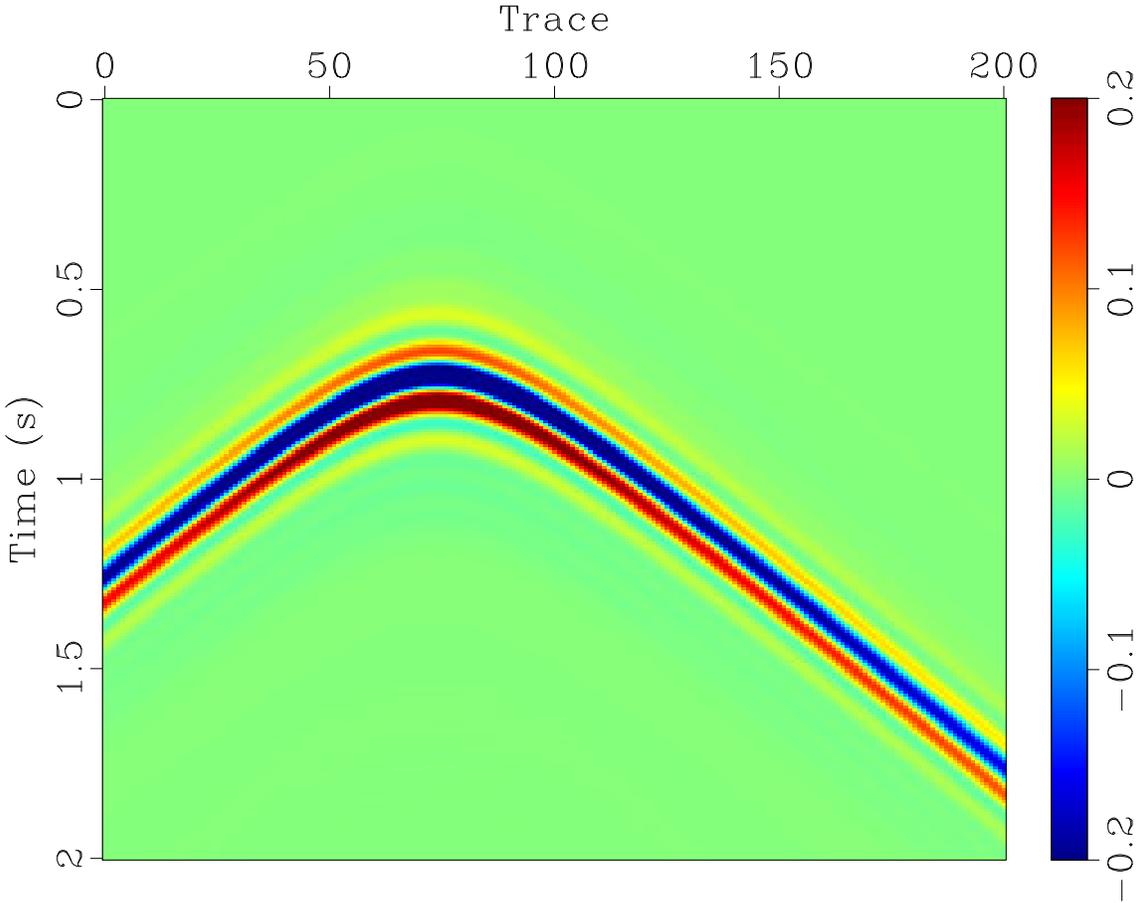
\plot{dfshh0}{width=\textwidth}{Velocity source gather = -2 $\times$
  pressure gather (Figure \ref{fig:dsrcphh0}).}

Figure \ref{fig:dinvhslh0} shows the approximate inversion (via the
first equation in display \ref{eqn:approxinv}) of the
pressure gather shown in Figure \ref{fig:drecplh0}, inverted in the homogeneous model
(rather than in the lens model used to generate the data). The result
differs greatly from the pressure source shown in Figure
\ref{fig:dhshh0}, as it must since it results from inversion in the
wrong material model.
Some dip filter effect is unavoidable and is caused by the aperture
limitation of the acquisition geometry: the steeper dips in the source
gather (Figure \ref{fig:dhshh0}) do not contribute to the data, nor to
the inversion. Also, the limited receiver aperture causes truncation
artifacts in the inversion. However, this result is an
accurate inversion: re-simulation (application of $S^{+}_{z_s,z_r}$) {\em using the same (homogeneous)
  model as used in the inversion} results in accurate recovery (Figure \ref{fig:drerecplh0}) of the
input pressure gather (Figure \ref{fig:drecplh0}). The difference is
shown on the same color scale in Figure \ref{fig:ddiffrecplh0}.

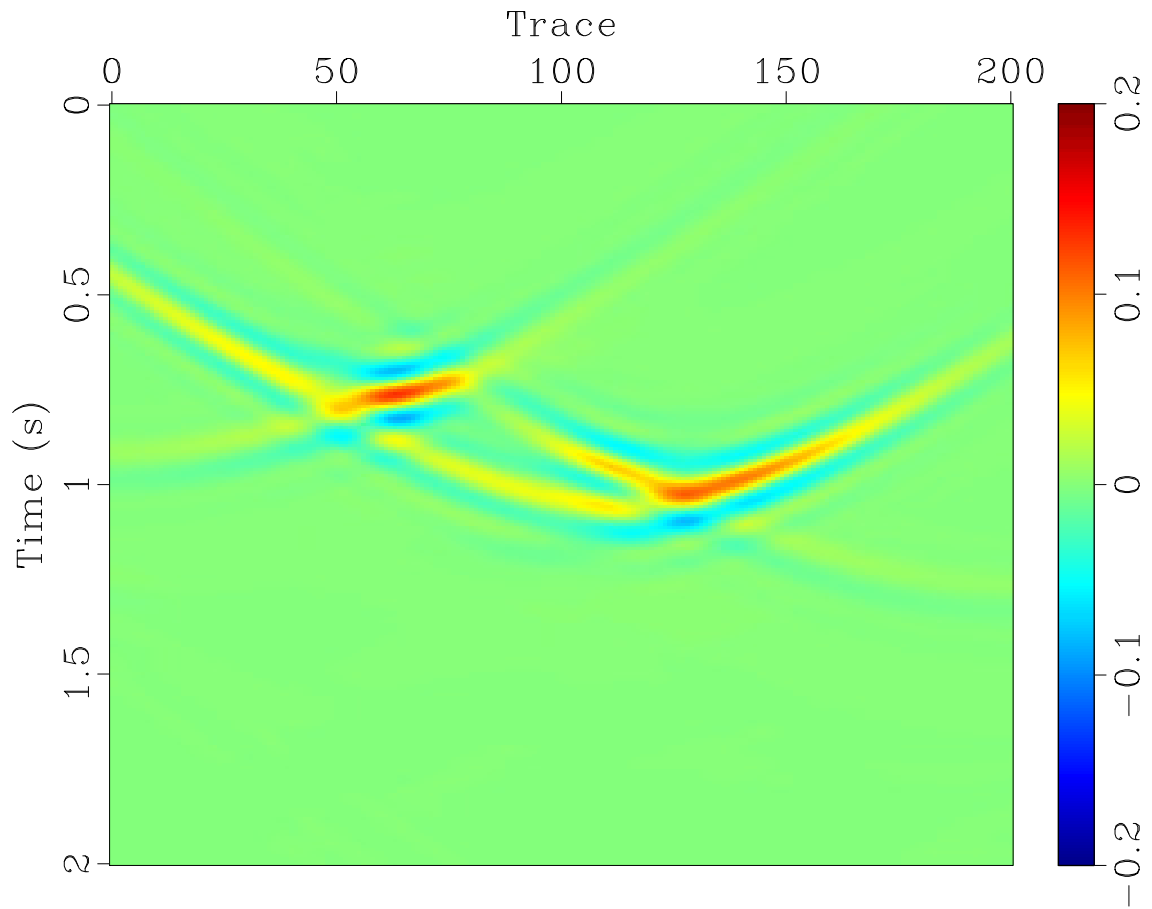
\plot{dinvhslh0}{width=\textwidth}{Approximate inversion via first 
  equation in display \ref{eqn:approxinv}. Inversion in homogenous model of 
  the pressure gather in Figure \ref{fig:drecplh0}, simulated with 
  lens model. Scaled version of output $v_z$ field obtained by 
  applying transpose of $4V^+_{z_s,z_r}$. Quite different from 
  pressure source gather (Figure \ref{fig:dhshh0}) used to generate
  data - since inversion takes place in a different material model!}

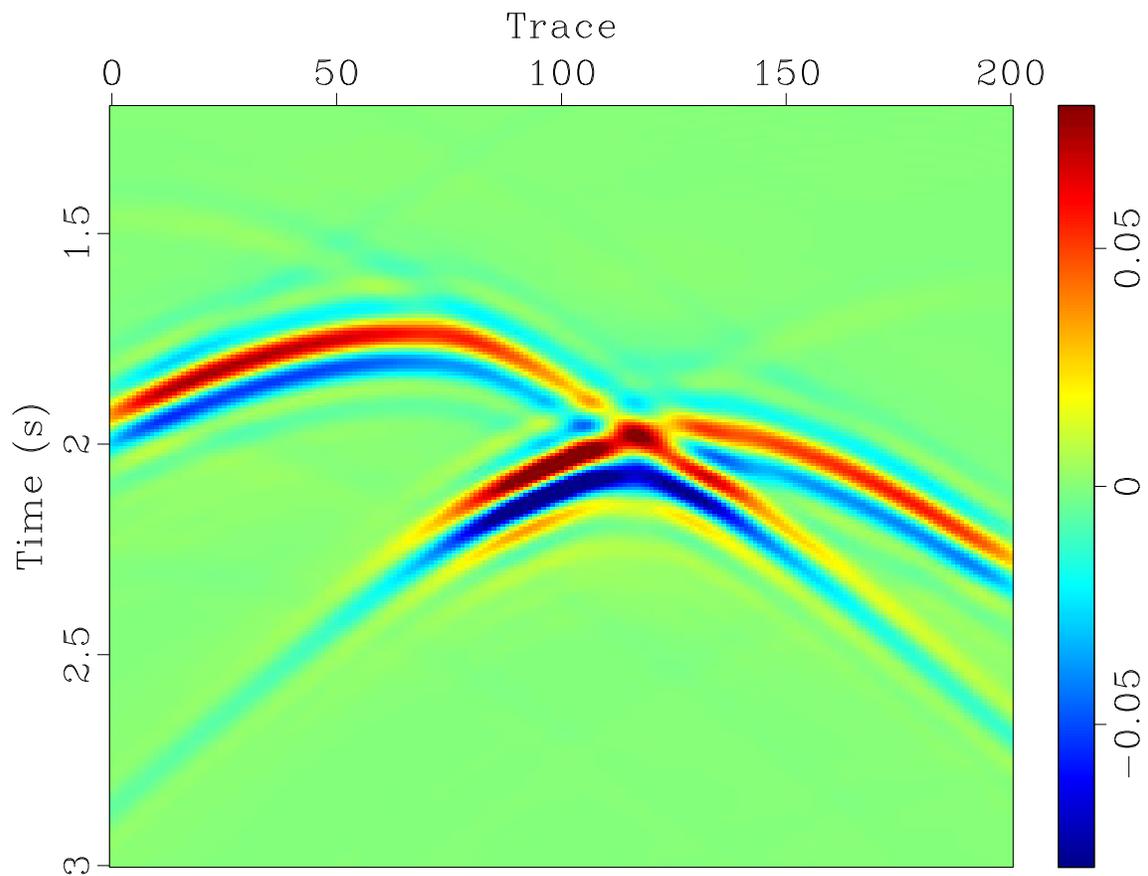
\plot{drerecplh0}{width=\textwidth}{Re-simulated pressure gather produced from
  inverted source
  shown in Figure \ref{fig:dinvhslh0}. Simulation in homogenous model
  used for inversion.}

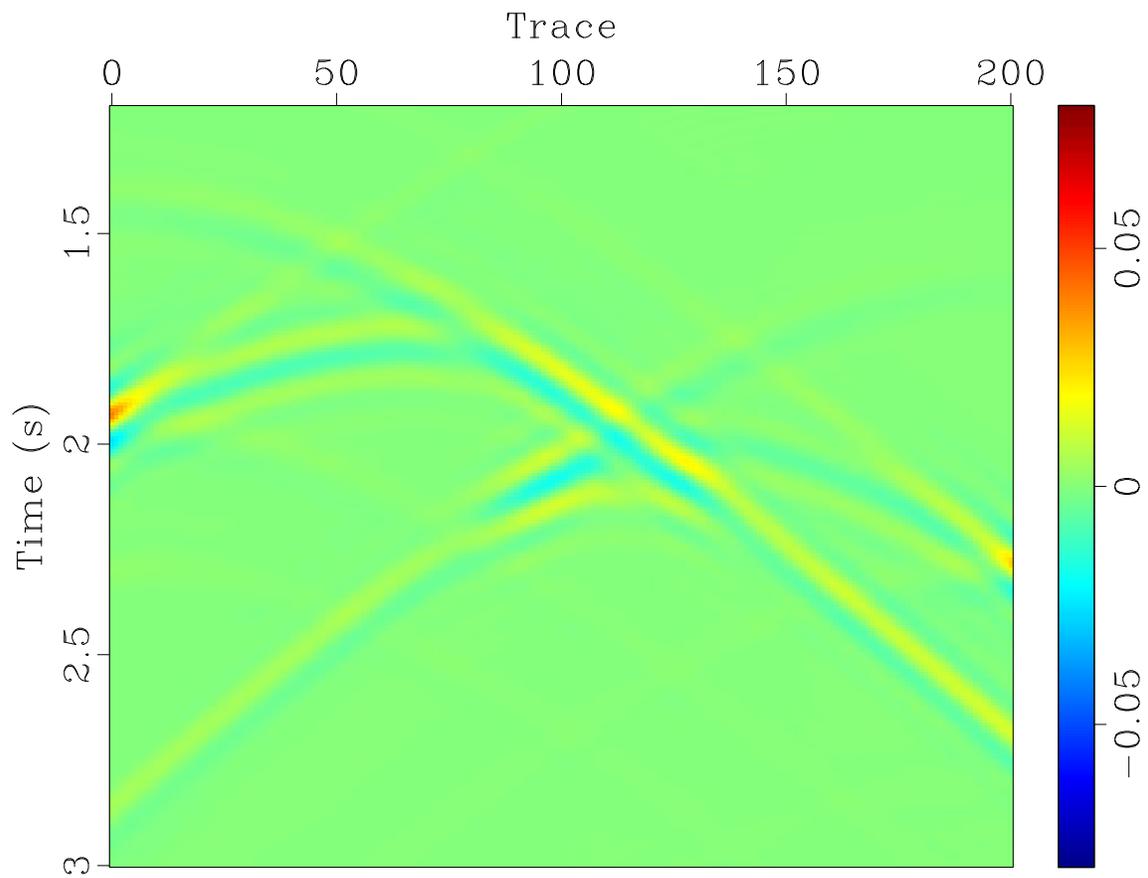
\plot{ddiffrecplh0}{width=\textwidth}{Difference between gathers
  displayed in Figures \ref{fig:drecplh0} and \ref{fig:drerecplh0},
  plotted on same color scale.}

\subsection{Quasi-unitary property of the modeling operator}
The identities \ref{eqn:approxinv} and \ref{eqn:sv} would together
establish the approximately unitary property of $S^+_{z_s,z_r}$, if
$\Lambda$ were symmetric. Identity \ref{eqn:approxinv} was illustrated
in the last subsection. Setting the symmetric issue aside for the
moment, an illustration of the relation \ref{eqn:sv} proceeds as follows.

Relation \ref{eqn:tracejump10} characterizes
$\Lambda^+_{z_s}$ as connecting the pressure and velocity components
of downgoing fields restricted to $z=z_s$. That is, the pressure
source gather
$h_s=\Lambda^{+}_{z_s}P_s p^+ = -2
P_{z_s}v_z$, displayed in Figure \ref{fig:dhshh0}, is the image of the
pressure gather in Figure \ref{fig:dsrcphh0} under
$\Lambda^{+}_{z_s}$.
The pressure gather
\ref{fig:dfwdplh0} is the image of this pressure source gather under
$S^+_{z_s,z_r}$ (using the lens model). The corresponding vertical velocity gather
(Figure \ref{fig:dfwdvzlh0}) is $-1/2$ times $\Lambda^+_{z_r}P_r
p$. Therefore scaling the data in Figure \ref{fig:dfwdvzlh0} by
$-2$ produces
$\Lambda^+_{z_r}S^+_{z_s,z_r}\Lambda^{+}_{z_s}P_s p^+$. On the other
hand, figure \ref{fig:daltvzlh0} shows the result of applying
$V^+_{z_s,z_r}$ to $f_s=-2P_s p^+$. Therefore scaling the gather in
Figure \ref{fig:daltvzlh0} by $-\frac{1}{2}$ produces $V^+_{z_s,z_r}P_zp^+$.
Since the data in Figures \ref{fig:dfwdvzlh0} and \ref{fig:daltvzlh0}
are essentially identical, the relation \ref{eqn:sv} holds for this
example.

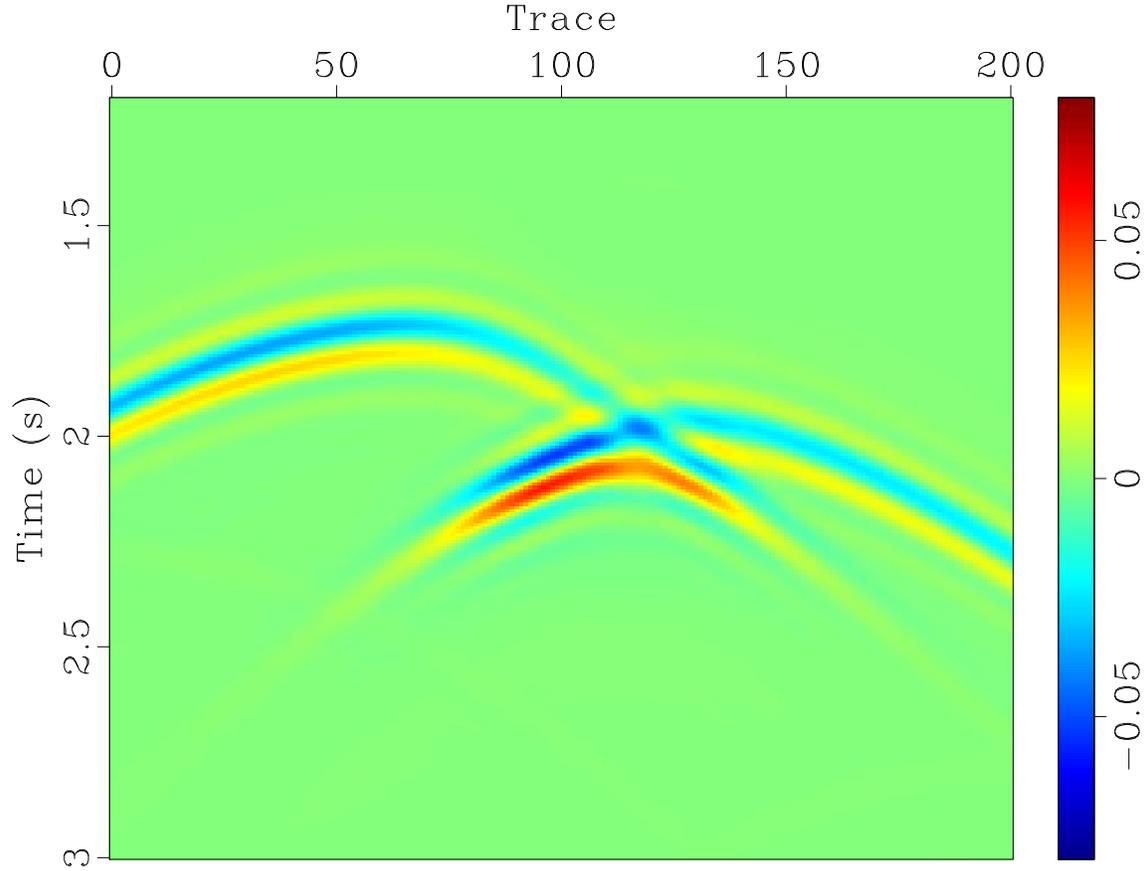
\plot{dfwdvzlh0}{width=\textwidth}{Vertical velocity gather, generated
  with a pressure source in the lens model, 
  corresponding to pressure gather \ref{fig:dfwdplh0}.}

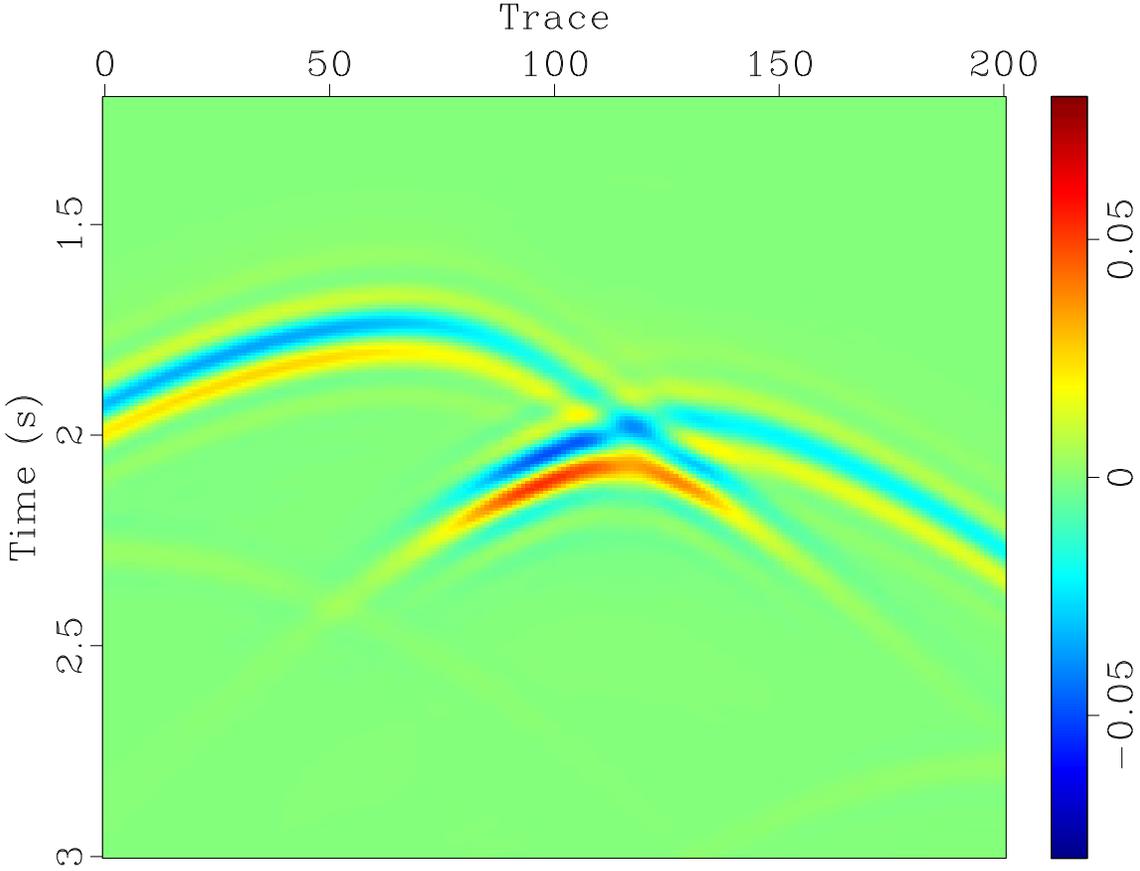
\plot{daltvzlh0}{width=\textwidth}{Vertical velocity gather, generated
  with a velocity source in the lens model, corresponding to
  pressure gather \ref{fig:daltplh0}.}

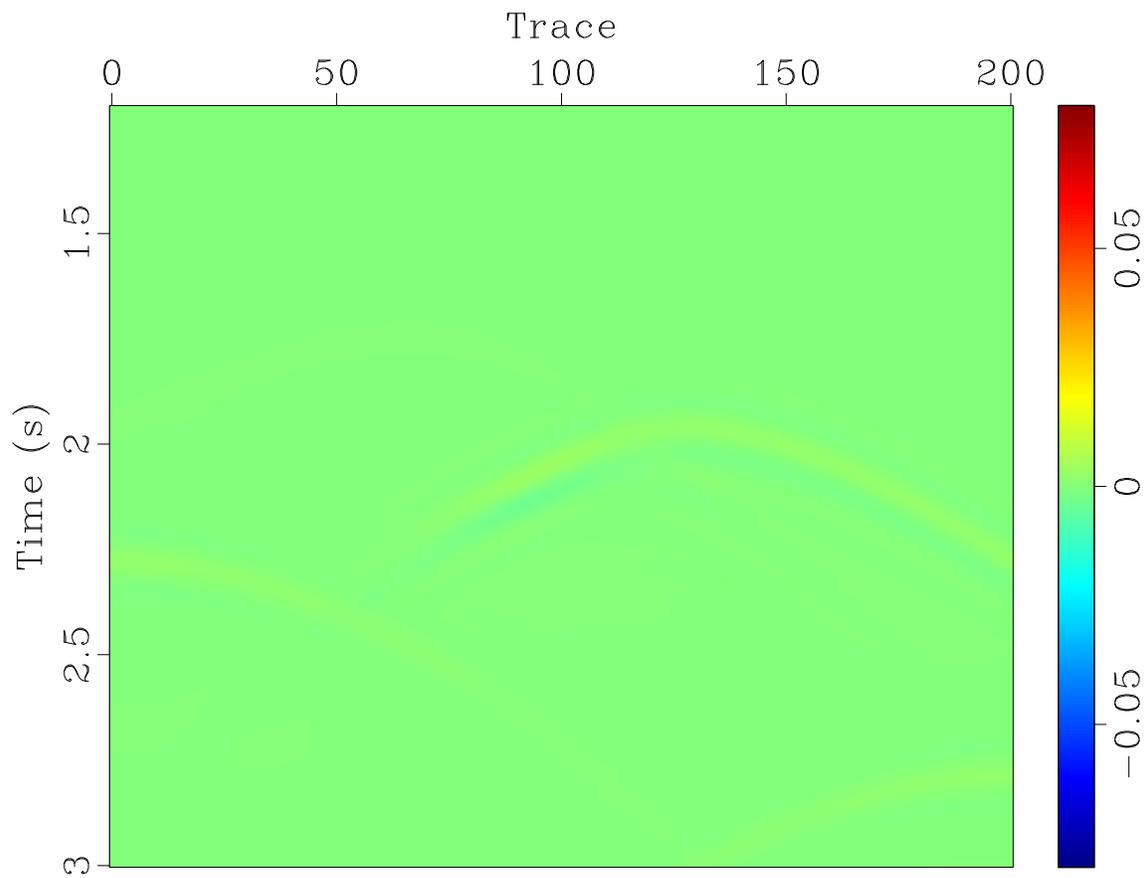
\plot{dsvcomplh0}{width=\textwidth}{Difference between velocity gathers
  shown in Figures \ref{fig:dfwdvzlh0} and \ref{fig:daltvzlh0},
  plotted on the same color scale as these figures.}



\subsection{Economical computation of $\Lambda$ in ``thin'' subdomain}
Equation \ref{eqn:lamnear} suggests a thin-slab computation of the
$\Lambda$ action, which is both accurate and economical. 
This calculations
place a receiver array at $z_s+\Delta z=2900$ m depth, just 100 m above the
source surface at $z_s=3000$ m. For the discretization used to create
the examples shown so far, that is just a 5 gridpoint difference in
depth, as opposed to 100 gridpoints between the source and receiver
depths for examples such as shown in Figure \ref{fig:drecplh0}.

Asymptotically, $\Lambda^{\pm}_{z_s}$ depends only on the medium
coefficients ${\bf c}$ in an arbirarily small region containing the
source surface $z=z_s$. In this example, the homogeneous and lens models are identical in the depth
range $2900 < z < 3000$ m, so the computed $\Lambda^+_{z_s}$ operators
will be precisely the same for both models. Hence I show only results for the
homogenous model.

The approximation to $\Lambda_{z_s}^+$ via equation \ref{eqn:lamnear}
for this configuration is evaluated in
Figures \ref{fig:preddnshshh0}, \ref{fig:ddiffnshshh0},
\ref{fig:dprednshsrecpll0}, and \ref{fig:ddiffprednshsrecpll0}. The effect of aperture
limitation is clearly diminished: the second figure in this series
compares the full-aperture pressure source gather (Figure
\ref{fig:dhshh0}) with the 
image of the corresponding pressure gather (Figure \ref{fig:dsrcphh0})
under the approximation to $\Lambda_{z_s}^+$, and the last figures
show that the approximated source gathers accurately predict the
point-source pressure gather at the receiver datum $z_r=1000$ m.

\plot{preddnshshh0}{width=\textwidth}{Pressure source gather = image
  under pressure-to-source operator $\Lambda^+_{z_s}$ of pressure gather
  shown in Figure \ref{fig:dsrcphh0}, homogeneous model, using
  ``near'' receiver traces at $z=2900$ m. Compare Figure
  \ref{fig:dhshh0}: because the sources and
  receivers are close, little aperture is lost in this case.}

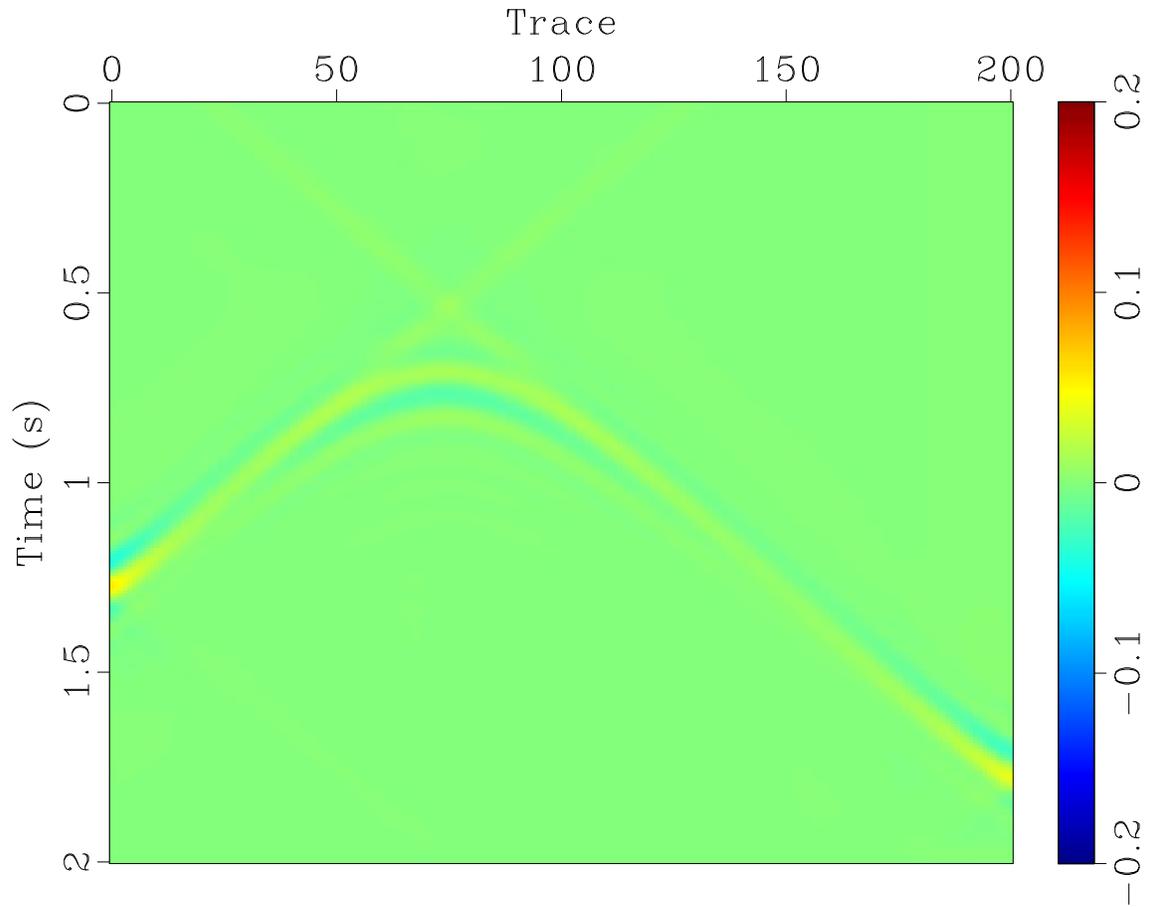
\plot{ddiffnshshh0}{width=\textwidth}{Difference between (a) image
  (Figure \ref{fig:preddnshshh0}) of $\Lambda^+_{z_s}$ applied to
  pressure gather (Figure \ref{fig:dsrcphh0}) using a near receiver
  array to implement formula \ref{eqn:lamident}, and (b) source gather
  (Figure \ref{fig:dhshh0}) inferred from vertical velocity.
  Homogeneous model used in all
  propagations. Same color scale as in Figure
  \ref{fig:preddnshshh0}. }



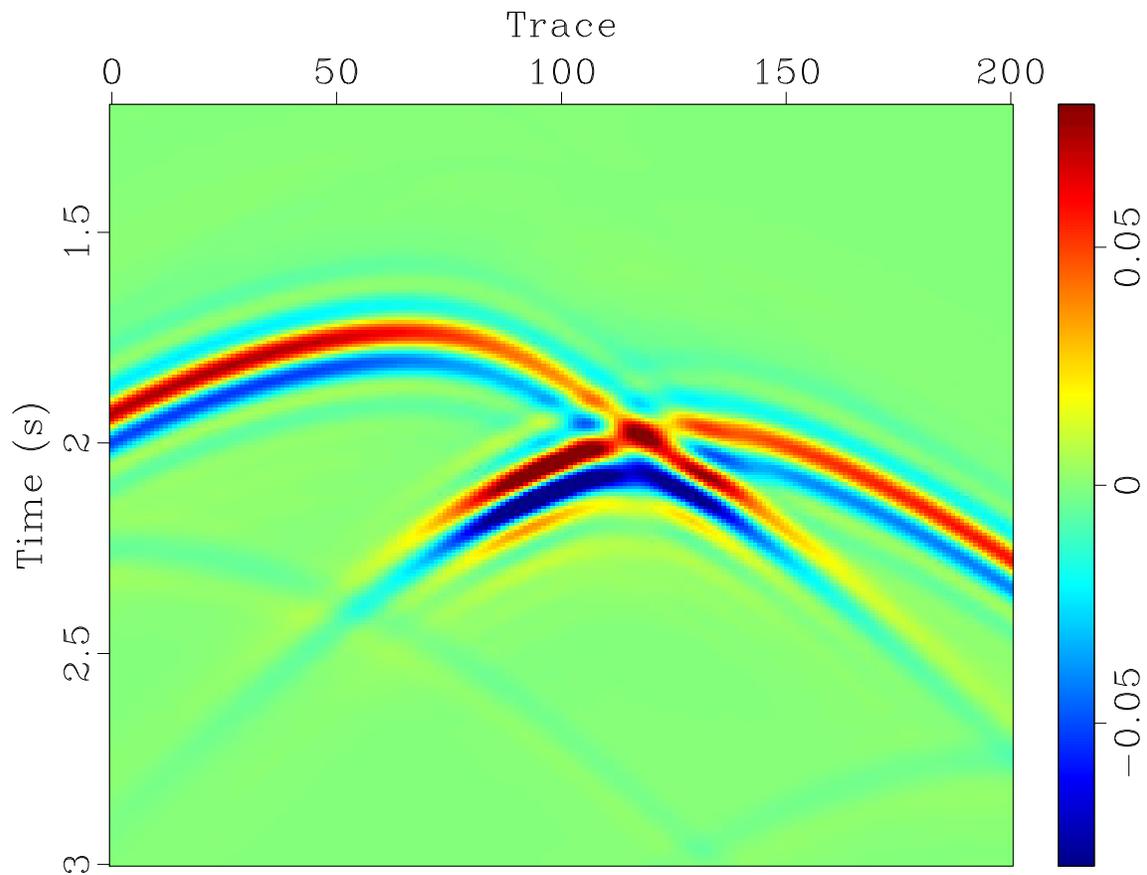
\plot{dprednshsrecpll0}{width=\textwidth}{Pressure gather at receiver
  datum $z=z_r=1000$ m simulated in lens model from source gather shown in Figure
  \ref{fig:preddnshshh0}. Compare with point source pressure gather
  (Figure \ref{fig:drecplh0}).}

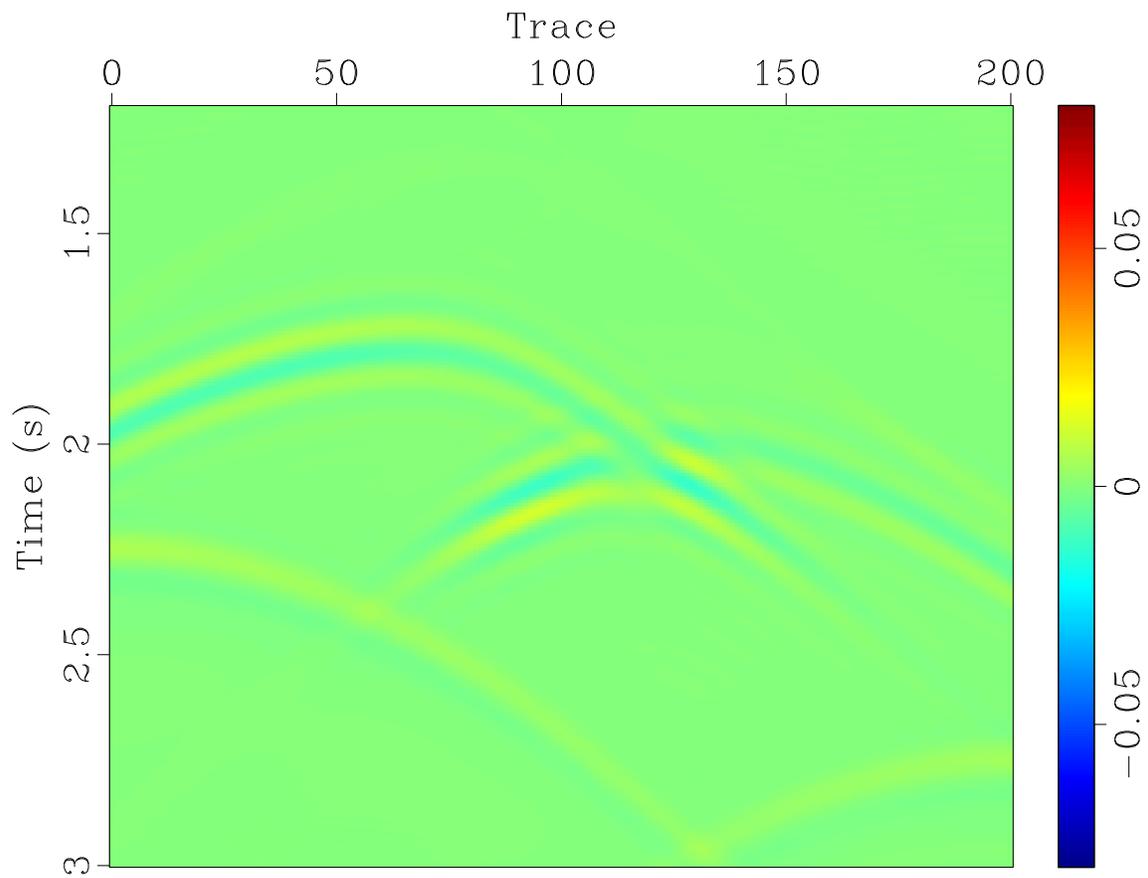
\plot{ddiffprednshsrecpll0}{width=\textwidth}{Plot of difference 
  between data shown in Figures \ref{fig:drecplh0} and 
  \ref{fig:dprednshsrecpll0}, plotted on same color scale as the latter 
  two figures.}
\subsection{Symmetrizing $\Lambda$}

Figure \ref{fig:preddnshstrhh0} shows the image of the pressure gather
in Figure \ref{fig:dsrcphh0} under $(\tilde{\Lambda}^+_{z_s})^T$,
using the ``near'' traces at $z=2900$, that is, $\Delta z = 100$ m in
expression \ref{eqn:lamtransp}, and propagation in the
homogeneous model. Note
the close resemblance to the image of the same pressure gather under
$\tilde{\Lambda}^+_{z_s}$ displayed in Figure
\ref{fig:preddnshshh0}. The difference of these two images is
displayed in \ref{fig:ddiffnslamtrhh0}, on the same color scale as the
images themselves. Since the propagation takes place entirely in a
region where all of the mechanical parameters are homogenous, I do not
offer a similar comparison for the lens model.

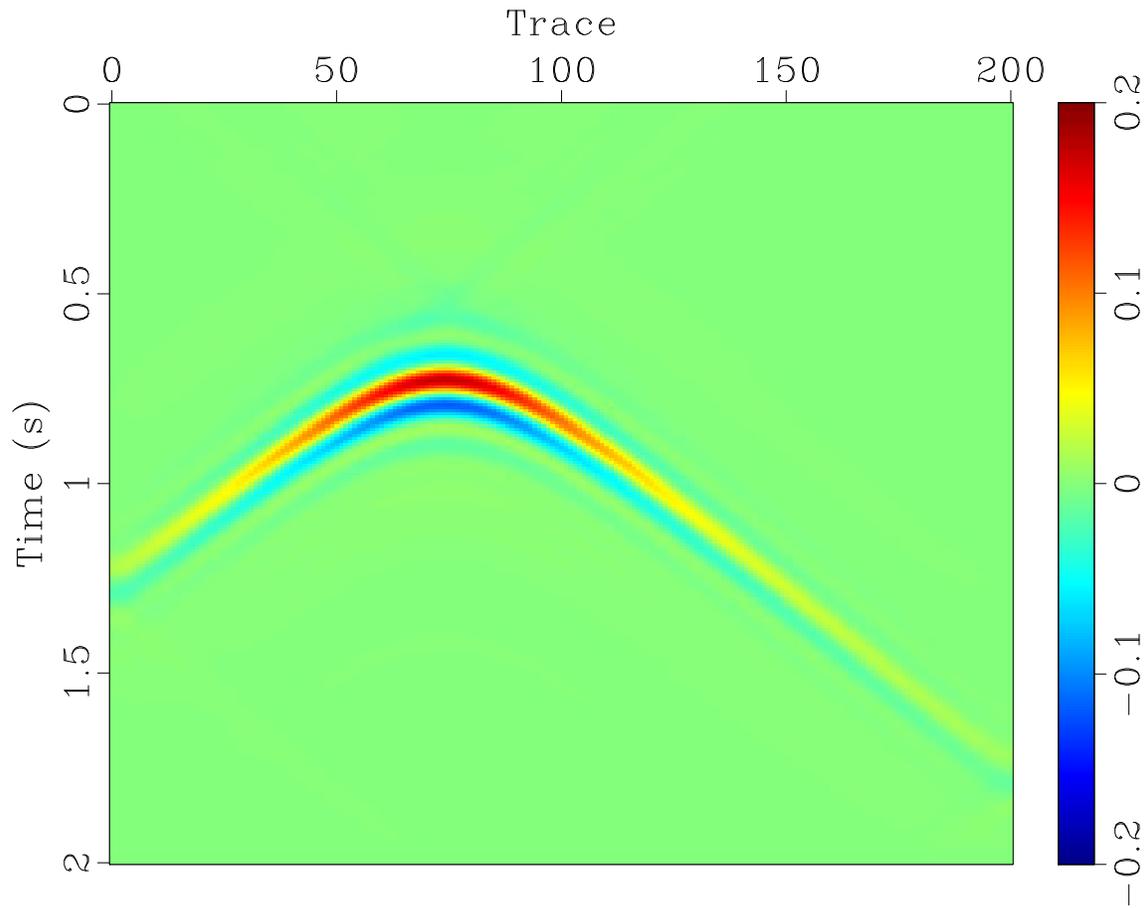
\plot{preddnshstrhh0}{width=\textwidth}{Pressure source gather = image
  under {\em transpose} of pressure-to-source operator
  $\Lambda^+_{z_s}$ of pressure gather shown in Figure
  \ref{fig:dsrcphh0}, homogeneous model, using ``near'' receiver
  traces at $z=2900$ m. Compare Figure \ref{fig:dhshh0} and
  \ref{fig:preddnshshh0}: as noted in the text, $\Lambda^+_{z_s}$ is
  asymptotically symmetric, so the resemblance is not a surprise.}

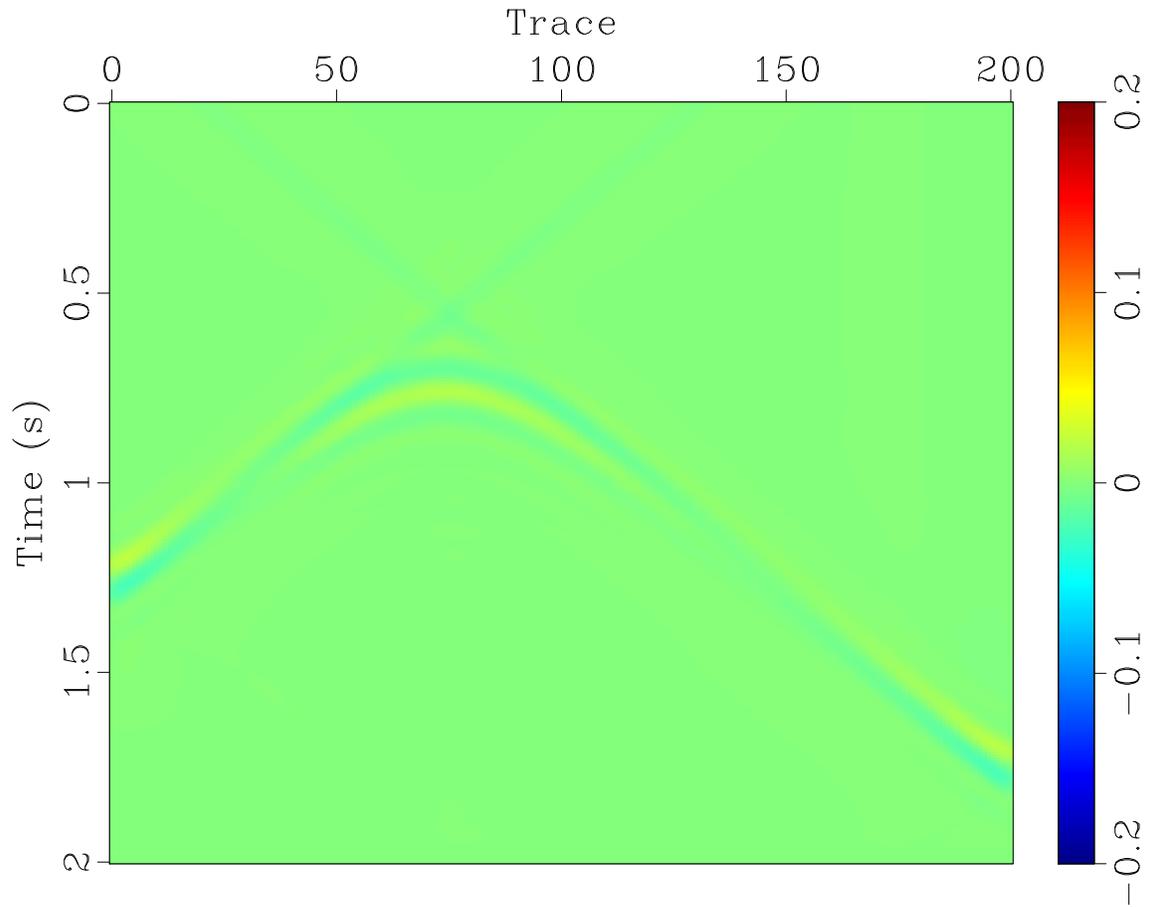
\plot{ddiffnslamtrhh0}{width=\textwidth}{Difference between data in
  Figures \ref{fig:preddnshshh0} and \ref{fig:preddnshstrhh0}, plotted
  on the same scale as these figures, showing that the asymptotic
  symmetry of $\Lambda^+_{z_s}$ is actually quantitative for the
  length, time and frequency scales of these examples.}

\subsection{Asymptotic symmetry of $\Lambda$}

Figure \ref{fig:symmdnshshh0} shows the output of the symmetrized
approximate source-to-pressure operator per equation \ref{eqn:winvcomp},
applied once again to the pressure data in Figure
\ref{fig:dsrcphh0}. Note the resemblance to Figures
\ref{fig:preddnshshh0} and \ref{fig:preddnshstrhh0}. These are all
asymptotic approximations of each other. Figure
\ref{fig:ddiffsymmdnshsll0} shows the
 the difference between the pressure gather at $z=z_r$ produced from
 the pressure source output by the symmetrized $\Lambda$, and the point source
simulation (Figure \ref{fig:dfwdplh0}), plotted on the same scale as
the latter, in both cases with all propagations in the lens model.

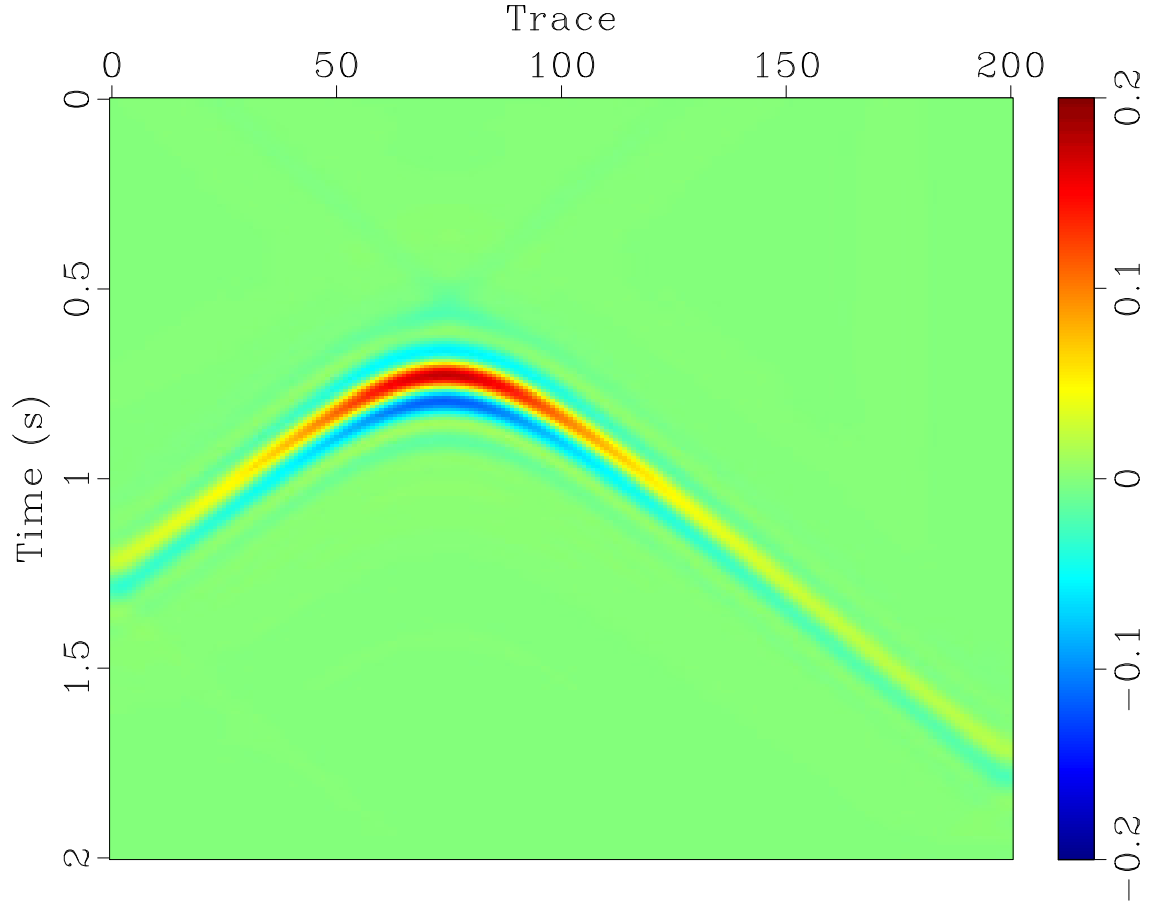
\plot{symmdnshshh0}{width=\textwidth}{Pressure source gather = image
    under {\em symmetrized} pressure-to-source operator
    $\frac{1}{2}\left(\Lambda^+_{z_s}+(\Lambda^+_{z_s})^T\right)$ of
    pressure gather shown in Figure \ref{fig:dsrcphh0}, homogeneous
    model, using ``near'' receiver traces at $z=2900$ m. Compare
    Figure \ref{fig:preddnshshh0}.}

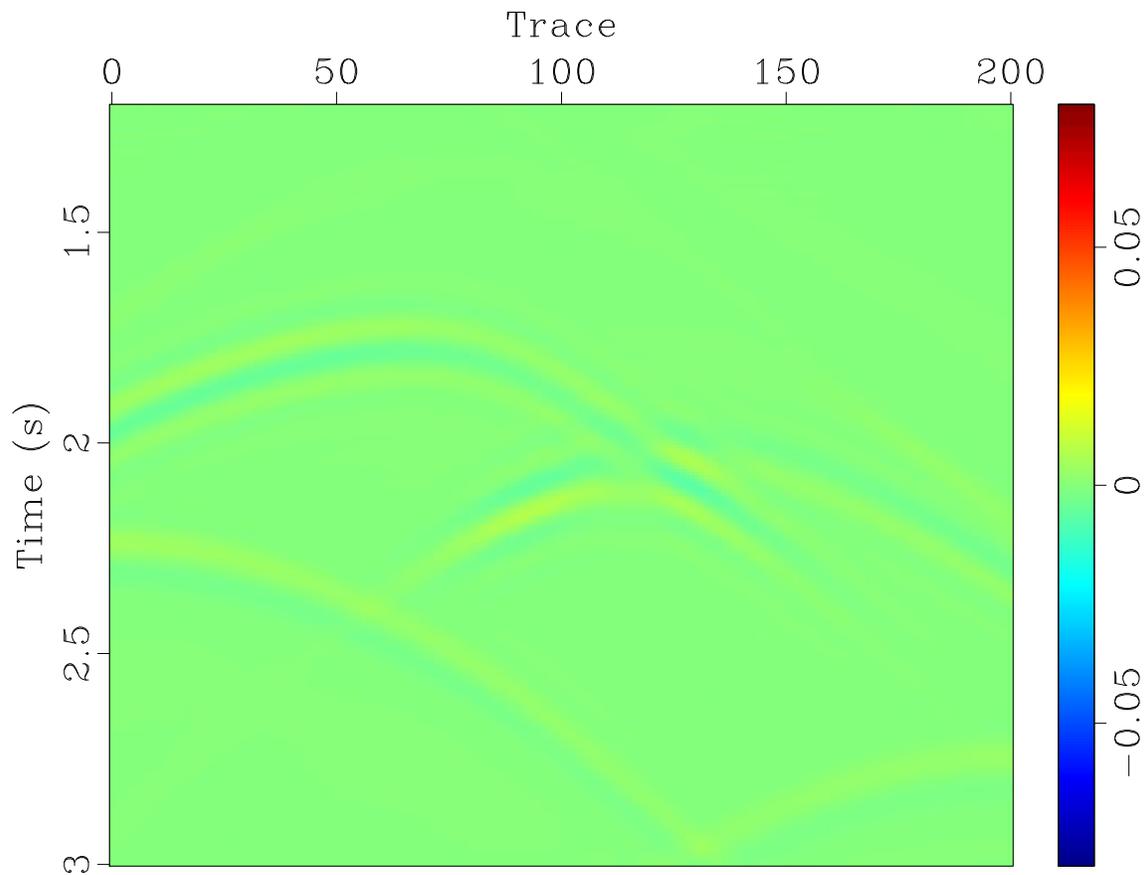
\plot{ddiffsymmdnshsll0}{width=\textwidth}{Difference between point
  source simulation (Figure \ref{fig:dfwdplh0}) and pressure gather at
  $z=z_r=1000$ m produced by simulation with the source shown in
  Figure \ref{fig:symmdnshshh0},
  propagation in the lens model.}

\subsection{Unitary property of modeling operator}
To illustrate this unitary property of $S^+_{z_s,z_r}$, I apply the
operator
\[
  \frac{1}{2}((\Lambda^+_{z_s})^T+ \Lambda^+_{z_s})
  (S^{+}_{z_s,z_r})^T \frac{1}{2}((\Lambda^+_{z_r})^T+
  \Lambda^+_{z_r})
\]
to the data $S^+_{z_s,z_r}h_s$ (Figure \ref{fig:drerecplh0}), in which
$h_s$ is the 
downgoing source created earlier (Figure \ref{fig:dhshh0})
The operator above is computed via the technique explained in
the preceding subsection, below, using
auxiliary receiver arrays 100 m above the data source and receiver arrays.

The output is shown in Figure
\ref{fig:lamsstlamrdrecplh0}. The difference with the actual source is
shown in Figure \ref{fig:difflamsstlamrdrecplh0}.

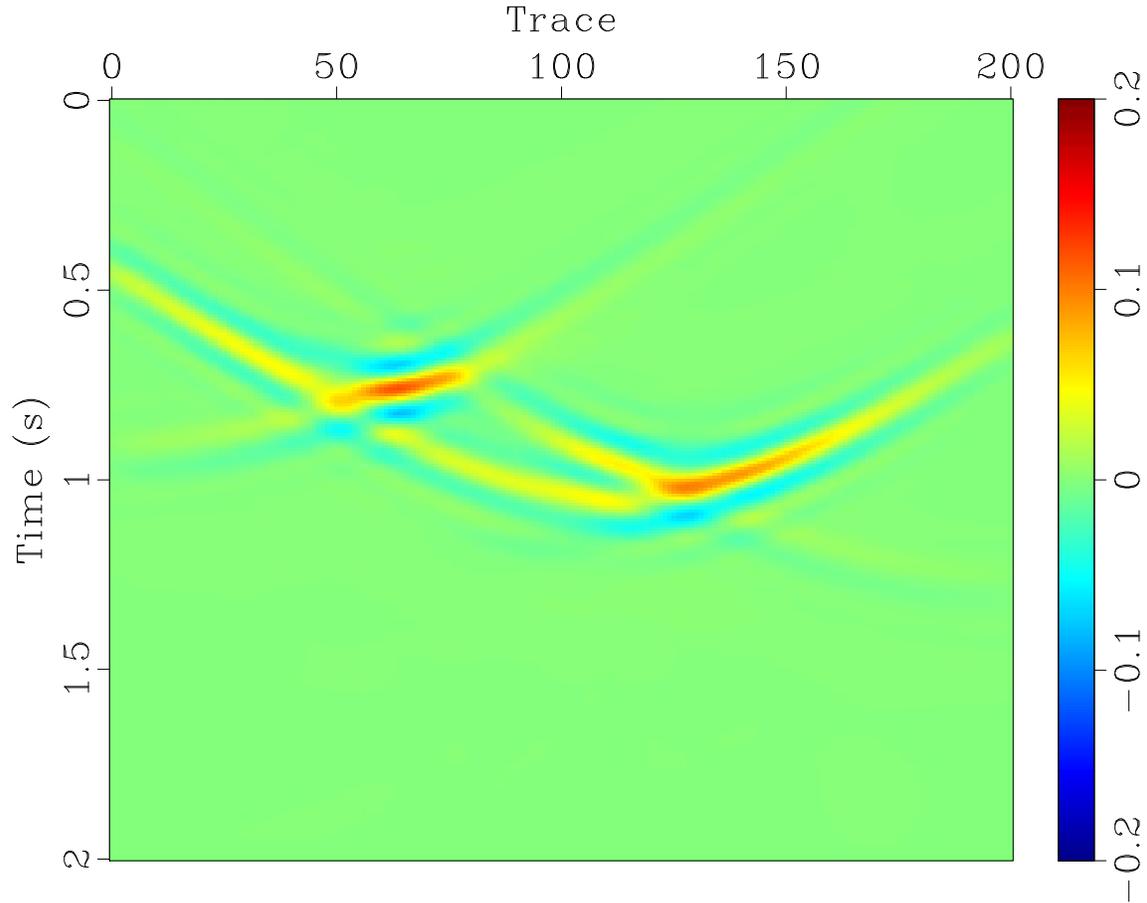
\plot{lamsstlamrdrecplh0}{width=\textwidth}{Inversion of data shown in
  Figure \ref{fig:drecplh0}, simulated in lens model, using the
  approximate unitarity relation \ref{eqn:unitary} and propagation in
  homogenous model.}

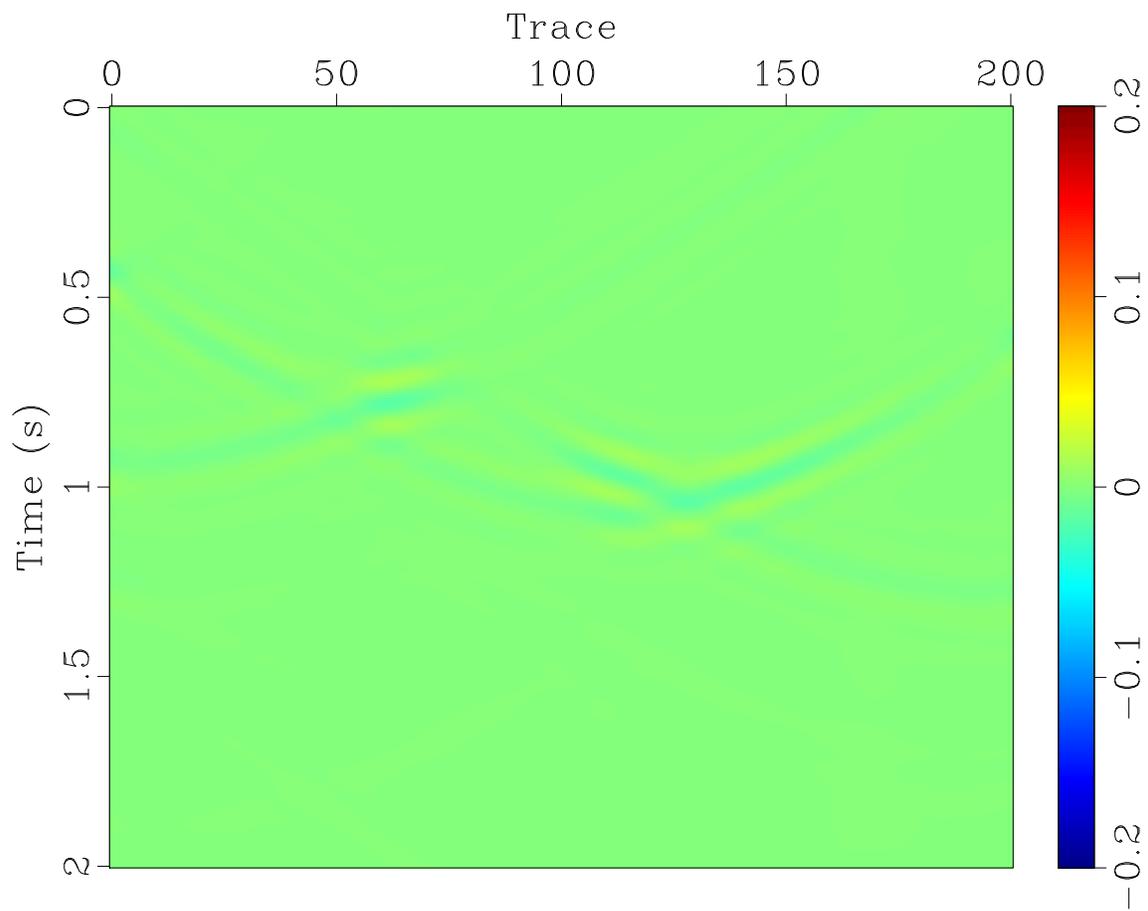
\plot{difflamsstlamrdrecplh0}{width=\textwidth}{Difference between
  data displayed in Figures \ref{fig:dinvhslh0} and
  \ref{fig:lamsstlamrdrecplh0}, plotted on the same color scale.}

\subsection{Preconditioned CG iteration}
This final subsection shows that result of Conjugate Gradient
iteration, with and without preconditioning, applied to the source
estimation problem \ref{eqn:esis}, with zero and non-zero penalty
weight $\alpha$. The data $d$ is the gather shown in
\ref{fig:drecplh0}, simulated using the lens model with source shown
in Figure \ref{fig:dhshh0}, or, alternatively, a point source with
bandpass filter wavelet located at $x_d=3500$ m, $z_d=3500$ m. In the inversion,
the material model is taken to be homogeneous, as has been the case in
all of the previous examples. 

Figure \ref{fig:compnres0lh0} shows the progress of the normal residual
(Euclidean norm of the difference of the two sides of equation \ref{eqn:norm1}),
for Conjugate Gradient and Preconditioned Conjugate Gradient
(Algorithm 1) iterations, applied to solution of the optimization
problem \ref{eqn:einv} with $\alpha=0$. For CG, the norms are both the ordinary
Euclidean norm, $W_m=W_d=I$. For PCG, $W_m$ and $W_d$ are given in
display \ref{eqn:wdef}, with the symmetrized $\Lambda$s computed as
indicated in the preceding subsections. Convergence for the
preconditioned algorithm is roughly 4 times as fast.

\plot{compnres0lh0}{width=\textwidth}{Comparison of normal residual 
  (gradient) Euclidean norms: CG (blue), PCG (red), plotted 
  vs. iteration. Data = lens model, point source (Figure 
  \ref{fig:drecplh0}), inversion in homogenous model. Penalty weight 
  $\alpha=0$.}

Figure \ref{fig:compnres1lh0} shows the same comparison with non-zero
penalty weight, $\alpha=10^{-3}$. The PCG normal residual curve is
almost identical with that in the $\alpha=0$ case, wheras the CG
convergence has slowed down noticeably, being about five times as slow
as the preconditioned algorithm.

\plot{compnres1lh0}{width=\textwidth}{Comparison of normal residual 
  (gradient) Euclidean norms: CG (blue), PCG (red), plotted 
  vs. iteration. Data = lens model, point source (Figure 
  \ref{fig:drecplh0}), inversion in homogenous model. Penalty weight 
  $\alpha=10^{-3}$.}

\section{Conclusion}

The linear modeling operator of surface source extended acoustic
waveform inversion is approximately invertible, and this paper has
shown how to approximately invert it. The construction is based on
reverse time propagation of data, as inspired by the literature on
photoacoustic tomography. However, since the input energy comes from a
surface source, rather than a pressure boundary value, the
pressure-to-source operator intervenes. It provides not just an
approximate inverse, but a definition of weighted norms in domain and
range spaces of the modeling operator, in terms of which that operator
is approximately unitary. Accordingly, Krylov space iteration defined
in terms of these weighted norms, or equivalently preconditioned
Conjugate Gradient iteration, gives a rapidly convergent solution
method for the linear subproblem.

The existence of an approximate unitary representation of the modeling
operator is not merely a computational convenience, however. It
reveals fundamental aspects of the operator's structure that enable
an explanation for the mitigation of cycle-skipping, a feature of the
{\em nonlinear} extended inverse problem. This fact echoes earlier observations
concerned a reflected wave inverse problem, involving a modeling
operator with a similar approximate inverse
\cite[]{tenKroode:IPTA14,Symes:IPTA14}. Also, the approximate inverse
leads to a stable computation of the gradient of the nonlinear
objective function \ref{eqn:esi}, resolving a difficulty first noted
also for reflected wave inversion \cite[]{KerSy:94}.

The transmission inverse problem figuring most prominently in
contemporary applied seismology is surely the Full Waveform Inversion
(FWI) of diving wave data. This is essentially the same problem as the
one discussed in this paper, and can be formulated and treated the
same way, at least for acoustic wave physics. All of the topics
treated here are open for elastic wave physics - the analogue of the
pressure-to-source map would is the map from surface velocity field to
corresponding constitutive defect, analogous to the elastic
Dirichlet-to-Neumann map investigated by \cite{Rachele:00}.

The underlying tool in the ideas developed here is geometric optics
(or ray theory), without which the very concept of downgoing waves
would be meaningless. The physics of actual earth materials includes
material heterogeneity on all scales, which appears to leave little
room for the assumption of scale separation underlying geometric
optics. Moreover, earth materials are anelastic, with elastic wave
energy being converted to and from thermal excitation, pore fluid
motion, and so on. A truly satisfactory understanding of inverse wave
problems will eventually need to accommodate heterogenity and
anelasticity beyond the current capabilities of the ray-based theory.

\section{Declarations}
The author received no funding from any source in the performance of
the wo rk reported here, nor does he have any financial or
non-financial interests relevant to the content of this article. The
computations described were developed in the Madagascar reproducible
research framework (www.reproducibility.org). The software source is
available on request from the author.

\append{Adjoint Computation}

The adjoint of ${\cal S}^+_{z_s,z_r}$ can be computed by a variant of
the adjoint state method, in this case a by-product of the
conservation of energy. This calculation leads to equation
\ref{eqn:sadj1}, from which the other statements about adjoints made in the second
section of the paper follow.

Suppose that $p^-,\bv^-$ solve \ref{eqn:awepm}
with $(h_s,f_s\bf{e}_z)\delta(z-z_s)$ replaced by
$ (h_r,f_r\bf{e}_z) \delta(z-z_r)$. Then
\[
0 = 
\left(\int\, dx\,dy\,dz\, \frac{p^+ p^-}{\kappa} +  
\rho \bv^+ \cdot \bv^- \right)|_{t \rightarrow \infty}
-
\left(\int\, dx\,dy\,dz\, \frac{p^+ p^-}{\kappa} +  \rho \bv^+ \cdot \bv^- \right)|_{t \rightarrow -\infty}
\]
\[
= 
\int_{-\infty}^{\infty} \,dt\, \frac{d}{dt}\left(\int\, dx\,dy\,dz\, \frac{p^+ p^-}{\kappa} +  \rho \bv^+ \cdot \bv^- \right)
\]
\[
= 
\int_{-\infty}^{\infty} \,dt\, \left(\int\, dx\,dy\,dz\, \frac{1}{\kappa} \frac{\partial p^+}{\partial t} p^- +  p^+ \frac{1}{\kappa}\frac{\partial p^-}{\partial t} \right.
\]
\[
+
\left. \rho \frac{\partial \bv^+}{\partial t} \cdot \bv^- + \rho \bv^+ \cdot \frac{\partial \bv^-}{\partial t} \right)
\]
\[
= 
\int_{-\infty}^{\infty} \,dt\, \left(\int\, dx\,dy\,dz\, \left(- \nabla \cdot \bv^+ + 
 h_s \delta(z-z_s)\right) p^- + p^+ \left(- \nabla \cdot \bv^- + 
 h_r \delta(z-z_r)\right) \right.
\]
\[
+
\left.  (- \nabla p^++f_s\bf{e}_z) \cdot \bv^- + \bv^+ \cdot (-\nabla
  p^- + f_r \bf{e_z}) \right)
\]
\[
= 
\int_{-\infty}^{\infty}\,dt\, \left(\int\, dx\,dy\,dz\, \left(- \nabla \cdot \bv^+ + 
 h_s \delta(z-z_s)\right) p^- + p^+ \left(- \nabla \cdot \bv^- + 
 h_r \delta(z-z_r)\right) \right.
\]
\[
+
\left.  p^+ (\nabla \cdot \bv^-) + (\nabla \cdot \bv^+) p^- 
  +f_s \delta(z-z_s) v_z^- + v_z^+f_r \delta(z-z_r) \right)
\]
after integration by parts in the last two terms. Most of what is left cancels, leaving 
\[
0 = \int_{-\infty}^{\infty}\,dt\,dx\,dy\, (h_sP_sp^-+f_zP_sv_z^-) +
( h_rP_rp^++f_rP_rv_z^+) = \langle (h_s,f_s), {\cal S}^-(h_r,f_r) \rangle+ \langle (h_r,f_r), {\cal S}^+_{z_s,z_r}(h_s,f_s) \rangle
\]
whence \ref{eqn:sadj1} follows immediately.

\bibliographystyle{seg}
\bibliography{../../bib/masterref}

\begin{thebibliography}{}
\itemsep0pt

\bibitem[Bao and Symes, 1991]{BaoSy:91b}
Bao, G., and W. Symes,  1991, A trace theorem for solutions of linear partial
  differential equations: Mathematical Methods in the Applied Sciences, {\bf
  14}, 553--562.

\bibitem[Courant and Hilbert, 1962]{CourHil:62}
Courant, R., and D. Hilbert,  1962, Methods of mathematical physics, volume ii:
  Wiley-Interscience.

\bibitem[Dafni and Symes, 2018]{DafniSymes:SEG18b}
Dafni, R., and W. Symes,  2018, Accelerated acoustic least-squares inversion:
  88th Annual International Meeting, Expanded Abstracts, Society of Exploration
  Geophysicists, 4291--4295.

\bibitem[Golub and van Loan, 2012]{Golub:2012}
Golub, G.~H., and C.~F. van Loan,  2012, Matrix computations, 4th ed.: Johns
  Hopkins University Press.

\bibitem[Griewank, 2000]{Griewank:book}
Griewank, A.,  2000, Evaluating derivatives: Principles and techniques of
  algorithmic differentiation: Society for Industrial and Applied Mathematics
  (Frontiers in Applied Mathematics 19).

\bibitem[Hasco{\"e}t and Pascual, 2013]{TapenadeRef13}
Hasco{\"e}t, L., and V. Pascual,  2013, The {T}apenade {A}utomatic
  {D}ifferentiation tool: {P}rinciples, {M}odel, and {S}pecification: {ACM}
  {T}ransactions {O}n {M}athematical {S}oftware, {\bf 39}.

\bibitem[Hou and Symes, 2016a]{HouSymes:16}
Hou, J., and W. Symes,  2016a, Accelerating extended least-squares migration
  with weighted conjugate gradient iteration: Geophysics, {\bf 81}, no. 4,
  S165--S179.

\bibitem[Hou and Symes, 2016b]{HouSymes:EAGE16}
--------, 2016b, Accelerating least squares migration with weighted conjugate
  gradient iteration: 78th Annual International Conference and Exhibition,
  Expanded Abstract, European Association for Geoscientists and Engineers,
  P104.

\bibitem[Hu et~al., 2007]{Habashy:07}
Hu, W., A. Abubakar, and T. Habashy,  2007, Application of the nearly perfectly
  matched layer in acoustic wave modeling: Geophysics, {\bf 72}, SM169--SM176.

\bibitem[Huang et~al., 2019]{HuangNammourSymesDollizal:SEG19}
Huang, G., R. Nammour, W. Symes, and M. Dolliazal,  2019, Waveform inversion by
  source extension: 89th Annual International Meeting, Expanded Abstracts,
  Society of Exploration Geophysicists, 4761--4765.

\bibitem[Kern and Symes, 1994]{KerSy:94}
Kern, M., and W. Symes,  1994, Inversion of reflection seismograms by
  differential semblance analysis: {A}lgorithm structure and synthetic
  examples: Geophysical Prospecting, {\bf 99}, 565--614.

\bibitem[Lasiecka, 1986]{Lasi:86}
Lasiecka, I.,  1986, Sharp regularity results for mixed hyperbolic problems of
  second order: Springer Verlag, volume~{\bf 1223} {\it of} Springer Lecture
  notes in Mathematics.

\bibitem[Lasiecka et~al., 1986]{LasLionsTrig:86}
Lasiecka, I., J.-L. Lions, and R. Triggiani,  1986, Non-homogeneous boundary
  value problems for second order hyperbolic operators: Journal de
  Math\'{e}matiques Pures et Appliqu\'{e}es, {\bf 65}, 149--192.

\bibitem[Lasiecka and Trigianni, 1989]{Lasi:87}
Lasiecka, I., and R. Trigianni,  1989, Trace regularity of the solutions of the
  wave equation with homogeneous boundary conditions and compactly supported
  data: Journal of Mathematical Analysis and Applications, {\bf 141}, 49--71.

\bibitem[Lax, 2006]{Lax:PDENotes}
Lax, P.~D.,  2006, Hyperbolic partial differential equations ({C}ourant
  {L}ecture {N}otes): American Mathematical Society.

\bibitem[Payne, 1975]{Payn:75}
Payne, L.,  1975, Improperly posed problems in partial differential equations:
  Lecture Note~22, CBMS, Society for Industrial and Applied Mathematics,
  Philadelphia.

\bibitem[Rachele, 2000]{Rachele:00}
Rachele, L.,  2000, Boundary determination for an inverse problem in
  elastodynamics: Communications in Partial Differential Equations, {\bf 25},
  1951--1996.

\bibitem[Stefanov and Uhlmann, 2005]{StefUhl:05}
Stefanov, P., and G. Uhlmann,  2005, Stable determination of generic simple
  metrics from the hyperbolic {D}irichlet-to-{N}eumann map: International
  Mathematics Research Notices, {\bf 17}, 1047--1061.

\bibitem[Stefanov and Uhlmann, 2009]{StefanovUhlmannIP:09}
--------, 2009, Thermoacoustic tomography with variable sound speed: Inverse
  Problems, {\bf 25}, 075011.

\bibitem[Symes, 2007]{Symes:06a-pub}
Symes, W.,  2007, Reverse time migration with optimal checkpointing:
  Geophysics, {\bf 72}, SM213--222.

\bibitem[Symes, 2014]{Symes:IPTA14}
--------, 2014, Seismic inverse problems: recent developments in theory and
  practice: Inverse Problems - from Theory to Application, Proceedings,
  Institute of Physics, 2--5.

\bibitem[Symes et~al., 2020]{Symes:SEG20}
Symes, W., H. Chen, and S. Minkoff,  2020, Full waveform inversion by source
  extension: why it works: 90th Annual International Meeting, Expanded
  Abstracts, Society of Exploration Geophysicists, 765--769.

\bibitem[Symes and Payne, 1983]{Symes:83}
Symes, W., and L.~E. Payne,  1983, Trace theorem for solutions of the wave
  equation and the remote determination of acoustic sources: Mathematical
  Methods in the Applied Sciences, {\bf 5}, 131--152.

\bibitem[Symes et~al., 2011]{GeoPros:11}
Symes, W., D. Sun, and M. Enriquez,  2011, From modelling to inversion:
  designing a well-adapted simulator: Geophysical Prospecting, {\bf 59},
  814--833.
\newblock (DOI:10.1111/j.1365-2478.2011.00977.x).

\bibitem[Tang et~al., 2013]{TangXuZhang:13}
Tang, B., S. Xu, and Y. Zhang,  2013, 3{D} angle gathers with plane-wave
  reverse time migration: Geophysics, {\bf 78}, no. 2, S117--S123.

\bibitem[ten Kroode, 2012]{tenKroode:12}
ten Kroode, F.,  2012, A wave-equation-based {K}irchhoff operator: Inverse
  Problems,  115013:1--28.

\bibitem[ten Kroode, 2014]{tenKroode:IPTA14}
--------, 2014, A {L}ie group associated to seismic velocity estimation:
  Inverse Problems - from Theory to Application, Proceedings, Institute of
  Physics, 142--146.

\bibitem[Virieux, 1984]{Vir:84}
Virieux, J.,  1984, {S}{H}-wave propagation in heterogeneous media: Velocity
  stress finite-difference method: Geophysics, {\bf 49}, 1933--1957.

\bibitem[Virieux and Operto, 2009]{VirieuxOperto:09}
Virieux, J., and S. Operto,  2009, An overview of full waveform inversion in
  exploration geophysics: Geophysics, {\bf 74}, no. 6, WCC127--WCC152.

\bibitem[Xu et~al., 2012]{XuWang:2012}
Xu, S., D. Wang, F. Chen, G. Lambaré, and Y. Zhang,  2012, Inversion on
  reflected seismic wave: SEG Technical Program Expanded Abstracts, 1--7.

\bibitem[Xu et~al., 2011]{XuZhangTang:11}
Xu, S., Y. Zhang, and B. Tang,  2011, 3{D} angle gathers from reverse time
  migration: Geophysics, {\bf 76}, no. 2, S77--S92.

\bibitem[Zhang et~al., 2014]{YuZhang:14}
Zhang, Y., A. Ratcliffe, G. Roberts, and L. Duan,  2014, Amplitude-preserving
  reverse time migration: from reflectivity to velocity and impedance
  inversion: Geophysics, {\bf 79}, S271--S283.

\bibitem[Zhang and Sun, 2009]{Zhang:SEG09}
Zhang, Y., and J. Sun,  2009, Practical issues of reverse time migration: True
  amplitude gathers, noise removal and harmonic-source encoding: Beijing
  International Geophysical Conference and Exposition, Expanded Abstracts,
  Society of Exploration Geophysicists, 204--209.

\end{thebibliography}

\end{document}